\documentclass{amsart}
\usepackage{amsmath}
\usepackage{amscd}
\usepackage{amsfonts}
\usepackage{amssymb}
\usepackage{latexsym}
\usepackage{pifont}
\usepackage{array}
\usepackage{xypic}
\newcommand{\ncm}{\newcommand}
\ncm{\beq}{\begin{equation}}
\ncm{\eeq}{\end{equation}}

\newtheorem{thm}{Theorem}[section]
\newtheorem{pro}[thm]{Proposition}
\newtheorem{lem}[thm]{Lemma}

\newtheorem{thm&def}[thm]{Theorem \& Definition}

\theoremstyle{definition}
\newtheorem{defi}[thm]{Definition}
\newtheorem{exa}[thm]{Example}

\theoremstyle{remark}
\newtheorem{rmk}[thm]{Remark}

\numberwithin{equation}{section}

\def\Set{\mathsf{Set}}

\def\tw{\mathsf{tw}}

\def\M{\mathsf{M}}
\ncm{\Cgd}{\mathsf{Cgd}}
\ncm{\BGD}{\mathsf{Bgd}}
\ncm{\BgdMap}{\mathsf{BgdMap}}
\ncm{\Bim}{\mathbb{BIM}}
\ncm{\YD}{\mathcal{YD}}
\ncm{\Gal}{\mathsf{Gal}}
\ncm{\Fib}{\mathsf{Fib}}
\ncm{\DHR}{\mathsf{DHR}}

\ncm{\F}{\mathcal{F}}
\ncm{\G}{\mathcal{G}}

\ncm{\V}{\mathcal{V}}
\ncm{\cM}{\mathcal{M}}
\ncm{\W}{\mathcal{W}}
\ncm{\Z}{\mathcal{Z}}
\ncm{\bicomC}{^C\M^C}

\ncm{\asso}{\mathbf{a}}
\ncm{\luni}{\mathbf{l}}
\ncm{\runi}{\mathbf{r}}

\ncm{\End}{\operatorname{End}}
\ncm{\Hom}{\operatorname{Hom}}
\ncm{\BiEnd}{\operatorname{BiEnd}}

\newcommand{\ci}{\circ}

\def\o{\otimes}
\def\x{\times}

\ncm{\amalgo}[1]{\underset{\scriptscriptstyle #1}{\o}}
\ncm{\ractB}{\underset{\scriptscriptstyle B}{\ract}}
\ncm{\ractT}{\underset{\scriptscriptstyle T}{\ract}}
\ncm{\mash}{\Pisymbol{psy}{35}}
\ncm{\mashed}[1]{\underset{\scriptscriptstyle #1}{\Pisymbol{psy}{35}}}

\ncm{\oA}{\amalgo{A}}
\ncm{\oB}{\amalgo{B}}
\ncm{\oC}{\amalgo{C}}
\ncm{\oR}{\amalgo{R}}
\ncm{\oT}{\amalgo{T}}
\ncm{\oL}{\amalgo{L}}
\ncm{\oS}{\amalgo{S}}
\ncm{\oN}{\amalgo{N}}
\ncm{\oH}{\amalgo{H}}
\ncm{\oQ}{\amalgo{Q}}
\ncm{\shrp}{\,\sharp\,}
\ncm{\Co}{\,\square\,}
\ncm{\coC}{\,\square _C\,}
\ncm{\coD}{\,\square _D\,}
\ncm{\coH}{\,\square _H\,}
\ncm{\coCe}{\,\square _{C^e}\,}
\ncm{\ex}[1]{\underset{\scriptscriptstyle #1}{\x}}

\ncm{\rarr}[1]{\stackrel{#1}{\longrightarrow}}
\ncm{\larr}[1]{\stackrel{#1}{\longleftarrow}}

\ncm{\op}{\mathrm{op}}
\ncm{\coop}{\mathrm{coop}}
\ncm{\co}{\mathrm{co}}
\ncm{\Fi}{\varphi}
\ncm{\OR}{\overrightarrow}
\ncm{\OL}{\overleftarrow}
\ncm{\cohom}{\mathbf{cohom}}
\ncm{\xp}{x_{\langle 0 \rangle} \ract \eta (x_{\langle 1 \rangle})_{(1)}}

\def\oneB{_{(1)}}
\def\twoB{_{(2)}}

\def\oneL{_{[1]}}
\def\twoL{_{[2]}}

\ncm{\coa}[1]{^{\langle #1\rangle}}
\ncm{\coabar}[1]{^{\langle \overline{#1}\rangle}}

\ncm{\I}{\mathcal{I}}

\def\du1{\hat 1}

\def\ract{\triangleleft}
\def\lact{\triangleright}

\ncm{\under}{\mbox{\rm\_}\,}
\ncm{\Cnt}{\mathsf{C}}

\ncm{\ZZ}{\mathbb{Z}}

\ncm{\Oo}{\mathcal{O}}
\ncm{\Ha}{\mathcal{H}}
\ncm{\Ve}{\mathcal{V}}

\ncm{\sqrM}{$\sqrt{\text{Morita}}$}
\ncm{\ld}{\,.\,}
\ncm{\ud}{\,^.\,}
\ncm{\Mor}{\underset{k}{\sim}}
\ncm{\into}{\hookrightarrow}
\ncm{\coinv}[1]{{\co\text{-}#1}}
\ncm{\lZ}{\overrightarrow{\mathcal{Z}}}
\ncm{\rZ}{\overleftarrow{\mathcal{Z}}}
\ncm{\rbeta}{\overleftarrow{\beta}}
\ncm{\rgamma}{\overleftarrow{\gamma}}
\ncm{\lbeta}{\overrightarrow{\beta}}
\ncm{\lgamma}{\overrightarrow{\gamma}}
\ncm{\fgp}{\mathrm{fgp}}
\ncm{\adj}{\dashv}
\ncm{\Mref}{{\scriptscriptstyle M\text{-}\mathrm{ref}}}

\begin{document}

\title{Scalar extension of bicoalgebroids}
\author[I. B\'alint]{Imre B\'alint}
\address{Research Institute for Particle and Nuclear Physics, Budapest}
\email{balint@rmki.kfki.hu}
%\thanks{Supported by the Hungarian Scientific Research Fund, OTKA T-}

\begin{abstract}
After recalling the definition of a \emph{bicoalgebroid}, we define comodules and modules
over a bicoalgebroid. We construct the monoidal category of comodules,
and define Yetter--Drinfel'd modules over a bicoalgebroid. It is proved that the Yetter--Drinfel'd 
category is monoidal
and pre--braided just as in the case of bialgebroids, and is embedded into the one--sided 
center of the comodule category. We proceed to define Braided Cocommutative Coalgebras (BCC) over a 
bicoalgebroid, 
and dualize the scalar extension construction of \cite{Brz-Mil} and \cite{Kor-Bim}, originally applied to
bialgebras and bialgebroids, to bicoalgebroids. A few classical
examples of this construction are given. Identifying the comodule
category over a bicoalgebroid with the category of coalgebras of the
associated comonad, we obtain a comonadic (weakened) version of Schauenburg's
theorem.  Finally, we take a look at the scalar extension and braided
cocommutative coalgebras from a (co--)monadic point of view.
\end{abstract}

\maketitle

\section{Introduction}
Bicoalgebroids were introduced by Brzezi\'nski and Militaru in
\cite{Brz-Mil} as the structure that dualizes bialgebroids (in fact,
Takeuchi's $\times _R$--bialgebras) in the sense of reversing
arrows. This notion is not to be confused with the different kinds
of bialgebroid--duals that were later introduced in \cite{K-Sz}. It
would seem that the study of bicoalgebroids hasn't been taken up
vigorously since their inception; in our view, they merit attention
for at least two reasons. First, it is well established that a
bialgebroid may be thought of as a non--commutative analogue of the
\emph{algebra of functions} on a groupoid. It follows that a
bicoalgebroid, in turn, should be regarded as a non--commutative
analogue of the groupoid itself. This raises the hope that classical
constructions on groupoids may find their non--commutative
generalizations more easily in the context of bicoalgebroids.
Secondly, just as bialgebroids play a fundamental role in depth--two
extensions of algebras, it is expected that bicoalgebroids feature
prominently in extensions of coalgebras (from a different approach, in 
\cite{Lar-D2} Kadison constructs bialgebroids from depth 2 extensions of coalgebras). To
complete the picture, the dual Hopf--Galois theory of \cite{Schn1}
for extensions of coalgebras should generalize (from bialgebras)
to bicoalgebroids, giving a dual version of bialgebroid--Galois
theory. Further work in this latter direction is deferred to a
subsequent publication.

Central to this paper is the introduction of scalar extension for
bicoalgebroids. Incidentally, the construction that was shown in
\cite{Kor-Bim} to be a non--commutative version of scalar extension
was defined (for Hopf--algebras) in \cite{Brz-Mil} -- alongside with bicoalgebroids.

\section{Bicoalgebroids; comodules and modules}

Throughout, $k$ will be a field and the category $\M = \M_k$ of
$k$--modules will serve as our underlying category. The unadorned
$\otimes$ always means $\otimes _k$.

We use the ubiquitous Sweedler notation for coproducts and coactions. For a coalgebra
$\langle C, \Delta, \varepsilon \rangle$, the coproduct $\Delta: C \to C \o C$ on elements is denoted 
$\Delta (c) = c_{(1)} \o c_{(2)}$, with an implicit finite summation understood, 
i.e.~$c_{(1)} \o c_{(2)} = \sum_i {c_{(1)}}^i \o {c_{(2)}}^i$. Quite similarly, a right $C$--coaction 
$\rho_M: M \to M \o C$ will be denoted $\rho_M (m) = m_{[0]} \o m_{[1]}$ and a 
left $C$--coaction $\lambda_N: N \to C \o N$ will be denoted $\lambda_N (n) = 
n_{[-1]} \o n_{[0]}$.

The category of bicomodules
over a $k$--coalgebra $C$ is monoidal with monoidal unit $C$ and
monoidal product the cotensor product over $C$. This category will
be referred to as $\langle {\bicomC} , \coC , C \rangle$.

In fact, if $C$ is a coalgebra over a ring $R$ which is flat as an $R$--module,
then the category of $R$--flat $C$--bicomodules is monoidal with monoidal product 
the cotensor product over $C$, and monoidal unit $C$.

We shall also use the following standard notations throughout the paper. The co--opposite coalgebra of a
coalgebra $\langle C, \Delta, \varepsilon \rangle$ is $C_{cop} = 
\langle C, \Delta_{cop}, \varepsilon \rangle$, with the coproduct $\Delta_{cop} (c) = 
\tw_{C,C}\, \ci\, \Delta (c) = c_{(2)}\, \o \, c_{(1)}$. In analogy to the concept of enveloping algebra, the
co--enveloping coalgebra of $C$ is $C^e = \langle C \o C_{cop}, \tw_{23} \,\ci\, (\Delta \o \Delta_{cop}), 
\varepsilon \o \varepsilon \rangle$.

\smallskip
Following \cite{Brz-Mil}, we recall the following (somewhat lengthy)
\begin{defi}\label{definicio}
A \emph{left} bicoalgebroid $\langle H, \Delta, \varepsilon, \mu, \eta, \alpha, \beta, C \rangle$
consists of
\begin{itemize}
\item a $k$--coalgebra $\langle H, \Delta_H, \varepsilon_H \rangle$
\item two coalgebra maps $\alpha: H \to C$ and $\beta: H \to C_{cop}$,
such that $\alpha$ and $\beta$ 'cocommute', i.e.
$\alpha (h\oneB) \o \beta(h\twoB) = \alpha (h\twoB) \o
\beta(h\oneB)$. These maps furnish $H$ with a ($C \o C$)--bicomodule
structure, such that  $(H; \lambda_L,\lambda_R; \rho_L, \rho_R) \in
{^{C \o C} \M ^{C \o C}}$. The four $C$--coactions are:
\begin{align}
\lambda_L (h) &= \alpha(h\oneB) \o h\twoB\,, \;\; \rho_L (h) = h\twoB \o \beta(h\oneB)\nonumber \\
\lambda_R (h) &= \beta(h\twoB) \o h\oneB\,, \;\; \rho_R (h) = h\oneB \o \alpha(h\twoB)\nonumber
\end{align}
\item $C$--bicomodule maps $\mu_H: H \coC H \to H$ and $\eta_H: C
  \to H$ (multiplication \& unit) making $(H, \lambda_L, \rho_L)$ an algebra in ${\bicomC}$,
\end{itemize}
\smallskip
subject to the following axioms:
\smallskip
\begin{enumerate}
\item The multiplication map $\mu: H \coC H \to H$ satisfies:
\beq\label{brzmil1}
\sum_i \mu (g^i \o h^i_{(1)}) \o \alpha(h^i_{(2)}) = \mu (g^i_{(1)} \o h^i) \o \beta(g^i_{(2)}) 
\eeq
\item and it is comultiplicative:
\beq\label{brzmil2}
\Delta \circ \mu (\sum_i g^i \o h^i) = \sum_i \mu (g^i_{(1)} \o
h^i_{(1)}) \o \mu (g^i_{(2)} \o h^i_{(2)})
\eeq
\item Furthermore, the product is counital (note that this axiom seems
to be missing in Ref. \cite{Brz-Mil}):
\beq
\varepsilon (g) \varepsilon (h) = \varepsilon \circ \mu (g \o h)
\eeq
\item The unit map $\eta: C \to H$ satisfies the unit axiom:
\beq \mu \circ (\eta \Co H) \circ \lambda_L = H = \mu \circ (H \Co \eta)
\circ \rho_L \eeq
\item The unit map is compatible with the coalgebra structure in the
  following sense:
\beq\label{etaalphaeta} \Delta (\eta (c)) = \eta (c)_{(1)} \o \eta
(\alpha (\eta(c)_{(2)})) = \eta (c)_{(1)} \o \eta (\beta
(\eta(c)_{(2)})) \eeq \beq \varepsilon (\eta (c)) = \varepsilon (c)
\eeq
\end{enumerate}
\end{defi}

In \cite{Brz-Mil}, it is first proved that the condition \ref{brzmil1}
makes sense, i.e.~the two sides of the equation are well defined
maps. This, in turn, implies that \ref{brzmil2} makes sense, which
boils down to $(\mu \Co \mu) \ci \tw_{23} \ci (\Delta \Co \Delta)$ being
a well--defined map. The condition \ref{brzmil1} on the multiplication
map may be rephrased by saying that $\mu$ factorizes through the
\emph{cocenter} of the $C$--bicomodule $^CH \Co H^C$, where the two
coactions are $\lambda_R$ and $\rho_R$. We define the cocenter for
bicomodules as follows. 
\begin{defi}\label{cocentrum}
Let $M \in {\bicomC}$ a $C$--bicomodule. Define the map 
\begin{align}
\Phi: M \o C^* &\to M \nonumber \\
m \o \varphi &\mapsto m_{[0]}\, \varphi (m_{[1]}) - m_{[0]}\, \varphi
(m_{[-1]}) \nonumber
\end{align}
where $C^*$ denotes the $k$-dual of the coalgebra $C$. Then, the cocenter of $M$ is defined by the
cokernel map $\zeta: M \to \Z(M)$, where
$$\xymatrix{
M \o C^* \ar[r]^-{\Phi} & M \ar[r]^-{\zeta} & \Z(M) 
}$$
Introduce also the epi--mono factorization $\Phi: M \o C^* \stackrel{e}{\rightarrow} J_M \stackrel{i}{\rightarrow} M$
\end{defi}
The cocenter satisfies the following universal property. Let 
\beq\nonumber
W_M = \{m_{[0]} \o m_{[1]} - m_{[0]} \o m_{[-1]} \,\,|\, m\in
M \} \subset M \o C
\eeq
then for all $k$--module maps $f: M \to N$ which satisfy  
\beq\label{egyenlet}
(f \o C) (W_M) = 0
\eeq
i.e.~$f(m_{[0]}) \o m_{[1]} = f(m_{[0]}) \o m_{[-1]}$, there is a unique
$f': \Z(M) \to N$ such that $f = f' \ci \zeta$:
$$\xymatrix{
M \ar[r]^-{f} \ar[d]_-{\zeta}& N \\
\Z(M) \ar@{-->}[ur]_-{f'} &
}$$
Indeed, applying $(N \o \varphi)$ to \ref{egyenlet}, we find that $(f \o \varphi) (W_M) = 0$
for all $\varphi \in C^*$, i.e.~$f$ annihilates $J_M$.

\medskip
If the coalgebra $C$ is locally projective as a $k$--module (see
\cite{Brz-Wis}, {\bf 42.9}), then \\$(\zeta \o C) (W_M) = 0$. To see
this, note that for $C$ locally projective, $(\zeta \o C) (W_M) = 0$
if and only if $(Id \o \varphi) \ci (\zeta \o C) (W_M) = 0$ for all
$\varphi \in C^*$. This, however, holds by the definiton of $\zeta$.

Thus, for locally projective $C$, a $k$--module map $f: M \to N$ factorizes through $\zeta: M
\to \Z(M)$ \emph{if and only if} $(f \o C)\, (W_M) = 0$. Since, throughout this paper, we are working over a field, 
it is in fact unnecessary to explicitly assume local projectivity: modules over a field are always free, hence
they are projective. A projective module is also locally projective.

\medskip
We apply the above definition to the bicomodule $^CH \Co H^C$ afforded
by the coactions $\lambda_R$ and $\rho_R$. For reference, the
bicomodule structure is
\begin{align}
&(\lambda_R \Co H): \sum_i g^i \Co h^i \mapsto \sum_i \beta (g^i_{(2)}) \Co
g^i_{(1)} \Co h^i \label{cotsq1}\\
&(H \Co \rho_R): \sum_i g^i \Co h^i \mapsto \sum_i g^i \Co
h^i_{(1)} \Co \alpha(h^i_{(2)})\label{cotsq2}
\end{align}
We make the following
\begin{defi}
$H \boxtimes H$ is the cocenter of the bicomodule $^CH \Co H^C$,
$$\xymatrix{H \coC H \o C^* \ar[r]^-{\Phi_2} & H \coC H
  \ar[r]^-{\zeta_2}& H \boxtimes H} $$
where $ \Phi_2 (g^i \Co h^i \o \varphi) = g^i \Co h^i_{(1)}\, \varphi (\alpha(h^i_{(2)}))
- g^i_{(1)} \Co h^i\, \varphi (\beta(g^i_{(2)}))$
\end{defi}  
Using \ref{cotsq1} and \ref{cotsq2}, the multiplication map $\mu: H \coC H \to H$
factorizes through $H \boxtimes H$, i.e.
$$\xymatrix{
H \coC H \ar[r]^-{\mu} \ar[d]_-{\zeta} & H  \\
H \boxtimes H \ar@{-->}[ur]_-{f'}&   
}$$
precisely if $\sum_i \mu (g^i \o h^i_{(1)}) \o \alpha(h^i_{(2)}) = \mu
(g^i_{(1)} \o h^i) \o \beta(g^i_{(2)})$ (condition \ref{brzmil1})
holds. 
\medskip\\
This construction can be seen as dual to that of the Takeuchi
product $\times_R$. For a left bialgebroid $A$, the submodule $A \times_R A \into A \oR A$ is
the center of the $R$--bimodule $r\cdot (A \o A)\cdot r' = A t(r) \o A
s(r)$. It is well--known that there is no well--defined multiplication
on $A \oR A$, but $A \times_R A$ is a ring with component--wise
multiplication. The dual result is that even though comultiplication
is not well-defined on $H \coC H$, the factor $H \boxtimes H$ becomes a
well--defined coalgebra. This ensures that \ref{brzmil2} is
well-defined.  \\

The reader may easily convince herself that these axioms are dual to
those of a left bialgebroid $\langle A, \mu_A, \eta_A, \Delta_A,
\varepsilon_A, s, t, R \rangle$ in the sense of reversing arrows and
making the following substitutions:  $\langle A, \mu_A, \eta_A \rangle
\leftrightarrow \langle H, \Delta_H, \varepsilon_H \rangle$,
$\{ \Delta_A, \varepsilon_A \} \leftrightarrow \{ \mu_H, \eta_H \}$,
$\{s, t\} \leftrightarrow \{\alpha, \beta\}$, $R \leftrightarrow C$.

A \emph{right} bicoalgebroid is a $C$--bicomodule algebra with the coactions $\lambda_R$ and $\rho_R$, 
i.e.~we require $(H, \lambda_R, \rho_R)$ to be a monoid in the category of $C$--bicomodules.  
The axioms dualize those of a \emph{right} bialgebroid (cf.~the Lemma below).  
\medskip\\

We note here a result in line with the duality between bialgebroids and bicoalgebroids. It is 
well--known that the simplest right bialgebroid over a ring $R$ is the enveloping algebra 
$R^e = R \o R^{op}$ (its opposite is a left bialgebroid). The following, dual statement provides our 
first example of a bicoalgebroid:
\begin{lem}
The co--enveloping coalgebra $C^e = C \o C_{cop}$ is a right bicoalgebroid with the following 
structure maps. The source-- and target maps are given by
$$\alpha: C \o C_{cop} \to C,\; c \o \bar{c} \mapsto c\, \varepsilon(\bar{c})\;\;\text{and}\;\; 
\beta: C \o C_{cop} \to C_{cop},\; c \o \bar{c} \mapsto \varepsilon(c)\, \bar{c}$$ 
Multiplication is
$$ \mu^e: C \o C_{cop} \Co C \o C_{cop} \to C \o C_{cop},\; c \o \bar{c} \Co d \o \bar{d} \mapsto 
d \varepsilon(c)\, \varepsilon(\bar{d}) \o \bar{c}$$
and the unit map is $\Delta_{cop}$, $\eta^e: C \to C \o C_{cop},\; \eta^e(c) = c_{(2)} \o c_{(1)}$.   
\end{lem} 
\begin{proof}
$C^e$ has the (right--bicoalgebroid type) $C$--bicomodule structure coming from 
$\lambda_R = (\beta \o C^e)\ci \Delta_{cop}$ and $\rho_R = (C^e \o \alpha) \ci \Delta$. Explicity, the
coactions are
\begin{align}
\lambda_R (c \o \bar{c}) &= \bar{c}_{(1)} \o (c \o \bar{c}_{(2)}) \nonumber\\
\rho_R (c \o \bar{c}) &= (c_{(1)} \o \bar{c}) \o c_{(2)} \nonumber
\end{align}
Multiplication and unit are then seen to be (left and right) $C$--bicomodule maps with respect to the 
coactions $\lambda_R$ and $\rho_R$. Associativity and the unit property are easy calculations.  
$c \o \bar{c} \Co d \o \bar{d} \in C \o C_{cop} \coC C \o C_{cop}$ means 
\beq\label{cotenzor}
(c_{(1)} \o \bar{c}) \o c_{(2)} \o d \o \bar{d} = (c \o \bar{c}) \o \bar{d}_{(1)} \o (d \o \bar{d}_{(2)})
\eeq
We have to prove that multiplication factorizes through $C^e \boxtimes C^e$, i.e.~$\alpha ((c \o \bar{c})_{(1)}) \o 
\mu ((c \o \bar{c})_{(2)} \Co (d \o
\bar{d})) = \beta ((d \o \bar{d})_{(1)}) \o \mu ((c \o \bar{c}) \o (d
\o \bar{d})_{(2)})$. Inserting definitions, this reads
\beq\label{coTak}
c_{(1)} \o \mu^e ((c_{(2)} \o \bar{c}) \Co (d \o \bar{d})) = \bar{d}_{(2)} \o \mu^e ((c \o \bar{c}) 
\Co (d \o \bar{d}_{(1)})) 
\eeq  
By the definition of $\mu^e$, \ref{coTak} simplifies to
\beq
c \o d\, \varepsilon(\bar{d}) \o \bar{c} = \bar{d} \o d\, \varepsilon (c) \o \bar{c} 
\eeq
which is a consequence of \ref{cotenzor} (applying the counit map twice). It now makes sense to 
demand the compatibility of multiplication and comultiplication,
\begin{align}
&(\mu^e \o \mu^e)\ci \tw_{2,3} \ci (\Delta^e \o \Delta^e) ((c \o \bar{c})((d \o \bar{d})) = 
(c_{(1)} \o \bar{c}_{(2)})(d_{(1)} \o \bar{d}_{(2)}) \o \nonumber\\
&\o (c_{(2)} \o \bar{c}_{(2)})(d_{(2)} \o \bar{d}_{(1)}) = (d_{(1)} \varepsilon(c)\varepsilon(\bar{d}) \o 
\bar{c}_{(2)}) \o (d_{(2)} \o \bar{c}_{(1)}) = \nonumber\\
&= \Delta^e \ci \mu^e ((c \o \bar{c}) \o (d \o \bar{d})) \nonumber
\nonumber 
\end{align}
We skip the proof of the remaining compatibilities, all of them being trivial calculations.
\end{proof}
Just as in the dual case (where $R^{e, op} = R^{op} \o R$ is a left bialgebroid), we also have that 
$C^e_{cop} = C_{cop} \o C$ is a left bicoalgebroid. The proof is entirely similar.

\subsection{Comodules over a bicoalgebroid}
Based on experience with bialgebroids and dualization arguments, it
may be expected that a more categorical approach to bicoalgebroids
leads to the study of it's category of comodules. 

\medskip
It is a well--known fact that a coalgebra map $\gamma: D \to C$
induces a $C$--bicomodule structure on $D$ such that $D$ becomes a
comonoid in ${\bicomC}$. The category of $D$--comodules is then
naturally constructed as a subcategory of ${\bicomC}$. We
specialize this remark to the case of a (left--) bicoalgebroid $H$
over $C$.    
Consider the coalgebra map $\varphi = (\alpha \o \beta)\ci \Delta: H \to C \o C_{cop}$. The left and right 
$C^e$--coaction induced by $\varphi$ on $H$ are $\lambda = (\varphi \o H)\ci \Delta$ and 
$\rho = (H \o \varphi)\ci \Delta$. Inserting the definition of $\varphi$ and comparing with the notation of 
\ref{definicio}, 
\begin{align}
\lambda &= (C\o \rho_L^{op})\ci \lambda_L: H \to (C\o C_{cop}) \o H\\
\rho &= (\rho_R\o C)\ci \lambda_R^{op}: H \to H \o (C\o C_{cop}) 
\end{align}
(by $\lambda_R^{op}$, for example, we mean the right $C_{cop}$--coaction corresponding to $\lambda_R$ through
the isomorphism $^C\M \simeq \M^{C_{cop}}$). The image of $\Delta_H$ lies in $H \coCe H$, i.e.~we can introduce
the $C^e$--bicomodule map $\bar{\Delta}: H \to H \coCe H$ with
$$\xymatrix{
\Delta_H: H \ar[r]^-{\bar{\Delta}} & H \coCe H \ar[r]^-{\bar{\iota}_{H,H}} & H \o H
}$$
(here, $\bar{\iota}_{M,N}: M \coCe N \to M \o N$ is the equalizer defining the cotensor product over $C^e$).
On the other hand, the $C^e$--bicomodule map $\varphi$ is the composite
$$\xymatrix{
\varepsilon_H: H \ar[r]^-{\varphi} &
C^e \ar[r]^-{\varepsilon \o \varepsilon} & k
}$$
since $\varphi$ is a coalgebra map. It is then straightforward to show that 
\begin{lem}
$\langle H, \bar{\Delta}, \bar{\varepsilon} = \varphi \rangle$ is a comonoid in $^{C^e}\M^{C^e}$.
\end{lem}
\begin{proof}
Coassociativity is trivial, and the counit property reads
$$\xymatrix{
H \coCe H \ar[d]_-{\varphi\coCe H} & H\ar[l]_-{\bar{\Delta}} \ar[r]^-{\bar{\Delta}}\ar[d]_-{=} & 
H \coCe H \ar[d]^-{H \coCe \varphi}\\
C^e \coCe H \ar[r]_-{l_H} & H & H \coCe C^e \ar[l]^{r_H}  
}$$
which commutes, because $(\varepsilon \o \varepsilon)\ci (\varphi(h_{(1)})) h_{(2)} = \varepsilon_H(h_{(1)}) h_{(2)} = 
h$, etc. 
\end{proof}
We can now define the category of comodules over $H$.
\begin{defi}
A left $H$--comodule over a left bicoalgebroid $H$ is a pair $\langle M, \delta_M \rangle$, where 
$M \in {^{C^e}\M^{C^e}}$, and $\delta_M: M \to H \coCe M$ is a $C^e$--bicomodule map for which 
$$\xymatrix{
M \ar[r]^{\delta_M} \ar[d]_-{\delta_M} & H \coCe M \ar[d]^-{\bar{\Delta} \coCe M} & & 
M \ar[r]^{\delta_M} \ar[dr]_-{l_M}& H \coCe M \ar[d]^-{\varphi \coCe M} \\
H \coCe M \ar[r]_-{H \coCe \delta_M} &H \coCe H \coCe M & & & C^e \coCe M  
}$$
making $\delta_M$ a coassociative \& counital coaction. 

\medskip
The category of $H$--comodules $^H\M$ has objects the left $H$--comodules, and the arrows
$f: \langle M, \delta_M \rangle \to \langle N, \delta_N \rangle$ are the $C^e$--bicomodule maps 
$f: M \to N$ such that  
$$\xymatrix{
M \ar[r]^-{f} \ar[d]_-{\delta_M}& N\ar[d]^-{\delta_N}\\
H\coCe M \ar[r]_-{H \coCe f} & H \coCe N
}$$  
\end{defi} 
Summarizing, a left bicoalgebroid $H$ over $C$ is simultaneously a
monoid in the category ${^C\M^C}$ (with coactions $(H, \lambda_L, \rho_L)$) 
and a comonoid in the category $^{C^e}\M^{C^e}$ (with coactions $(H,
\lambda, \rho)$). This phenomenon is already familiar from the theory
of bialgebroids, namely that the algebra and coalgebra structures live in different monoidal categories.  

\medskip
The forgetful functor associated to the map $\varphi: H \to C^e$,
\begin{align}
&F: {^H\M} \to {^{C^e}\M} \simeq {^C\M^C}\\
&\langle M, \delta_M \rangle \to \langle M, (\varphi \o M)\ci \delta_M \rangle
\end{align}
is faithful and \emph{left} adjoint to $H \coCe \_\,: {^{C^e}\M} \to
{^H\M}$. Let us briefly recall the dual situation: 
a left bialgebroid $A$ over $R$ is an $R^e$--ring with $s \o t: R \o R^{op} \to A$, i.e.~a monoid in
$_{R^e}\M_{R^e}$. The forgetful functor $U: {_A\M} \to {_{R^e}\M}$ is \emph{right} adjoint to $A \o_{R^e}\,\_ : 
_{R^e}\M \to {_A\M}$. Furthermore, Schauenburg's theorem states that bialgebroid
structures on the $R^e$--ring $A$ are in one-to-one correspondance
with monoidal structures on the category $_A\M$ such that the
forgetful functor $U$ is strict monoidal. At this point, the question arises whether
a dual of this theorem holds for bicoalgebroids, namely: is there a one-to-one correspondance between
bicoalgebroid structures on the coalgebra $\langle H, \bar{\Delta}, \bar{\varepsilon} \rangle$ and monoidal 
structures on the category ${^H\M}$ such that $F: {^H\M} \to {^C\M^C}$ is strict monoidal? The next theorem gives 
the forward implication. We take up this question again in Section 4,
and look at the reverse implication from a comonadic point of view.    
\begin{thm}\label{Hcomod}
Let $H$ be a left bicoalgebroid over $C$. Then there is a monoidal
structure on $^H\M$ making the forgetful functor $F: {^H\M} \to
{^{C^e}\M} \simeq {^C\M^C}$ strict monoidal. Identifying $H$--comodules with their
underlying $C$--bicomodules, the monoidal product is $\coC$, the
cotensor product over $C$ and $C$ is the monoidal unit. 
\end{thm}
\begin{proof}
Assume there is a monoidal structure $\langle {^H\M}, \odot, I \rangle$ on ${^H\M}$ such that the forgetful
functor is strict monoidal, meaning that we have a triple $\langle F, F^2, F^0 \rangle$, where the maps 
$F^{M,N}: F (M\odot N) \to F(M) \coC F(N)$ and $F^0: F(I) \to C$ are identities. This is tantamount to specifying
\begin{itemize}
\item an $H$--comodule structure on $C$, 
$$ \delta_C: C \to H \coCe C, \;\; \text{and}$$ 
\item an $H$--comodule structure on the cotensor product of objects $M,N \in {\bicomC}$,
$$\delta_{M\Co N}: M \coC N \to H \coCe (M \coC N), $$
natural in $M$ and $N$
\end{itemize}
The bicoalgebroid structure on $H$ allows us to construct such maps $\delta_C$ and $\delta_{M\Co N}$. 

\medskip
The unit map $\eta: C \to H$ provides the desired $H$--comodule
structure on $C$: 
\begin{align}\nonumber
\delta_C = (H \o \alpha) \ci \Delta \ci \eta, \;\;  \delta_C (c) = \eta(c)_{(1)} \o \alpha (\eta(c)_{(2)})
\end{align}
This is indeed a coaction, 
\begin{align}
& (H \o \delta_C)\ci \delta_C (c) = 
\eta(c)_{(1)} \o \eta(\alpha(\eta(c_{(2)})))_{(1)} \o  \alpha
((\eta(\alpha(\eta(c_{(2)})))_{(1)})_{(2)}) = \nonumber \\
&= \eta (c)_{(1)} \o \eta (c)_{(2)(1)} \o \alpha (\eta(c)_{(2)(2)}) =
\eta (c)_{(1)(1)} \o \eta (c)_{(1)(2)} \o  \alpha (\eta (c)_{(2)}) =
\nonumber\\
&= (\Delta_H \o C) \ci \delta_C (c), \nonumber
\end{align}
applying \ref{etaalphaeta} in the second equality and coassociativity
in the third.

\medskip
For $M, N \in {^H\M}$, define the coaction $\delta_{M \Co N}: M \coC N \to H \coCe (M \coC N)$ as the 
composite map:
$$\xymatrix{
\delta_{M\Co N}:\; M \coC N \ar[r]^-{{\delta_M \coC \delta_N}}& (H \coCe M) \coC
(H \coCe N) \ar[r]^-{{\kappa}}& H \coCe (M \coC N)
}$$

Implicit in this definition is the map 
\beq
\kappa: (H \coCe M) \coC
(H \coCe N) \to  H \coCe (M \coC N)
\eeq
which we define as the unique arrow in the following diagram
$$\xymatrix{
(H \coCe M) \coC (H \coCe N) \ar@{-->}[rr]^-{\kappa} \ar[d]_-{\iota_{H\Co M, H\Co N}} & & H \coCe (M \coC N) 
\ar[d]^-{\bar{\iota}_{H, M\Co N}} \\
(H \coCe M) \o (H \coCe N) \ar[d]_-{\bar{\iota}_{H, M} \o \bar{\iota}_{H, N}} & & H \o (M \coC N) \ar[d]^-{H \o 
\iota_{M,N}}\\
(H \o M) \o (H \o N) \ar[rr]_-{(\mu_H \o M \o N)\,\ci \, \tw_{23}}  & & H \o (M \o N) 
}$$

By the definition of the kernel maps $\iota_{M,N}: M \coC N \to M \o
N$ and $\bar{\iota}_{U,V}: U \coCe V \to U \o V$, $(h \o m) \o (h' \o
n) \in {(H \coCe M) \coC (H \coCe N)}$ if and only if the following
identities hold:
\begin{align}
(h\oneB \o \alpha (h\twoB) \o m) \o (h' \o n) &= (h \o m_{[-1]} \o m_{[0]}) \o (h' \o n) \label{A}\\
(h\oneB \o \beta (h\twoB) \o m) \o (h' \o n) &= (h \o m_{[1]} \o m_{[0]}) \o (h' \o n) \label{B} \\
(h \o m) \o (h'\oneB \o \alpha (h'\twoB) \o n) &= (h \o m) \o (h' \o  n_{[-1]} \o n_{[0]}) \label{C}\\
(h \o m) \o (h'\oneB \o \beta (h'\twoB) \o n) &= (h \o m) \o (h' \o n_{[1]} \o n_{[0]}) \label{D}
\end{align}
and 
\beq
h\twoB \o m \o \beta (h\oneB) \o h' \o n = h \o m \o
  \alpha(h'\oneB) \o h'\twoB \o n \label{E}
\eeq

The arrow $\kappa$ is defined by the universal property of the composite kernel map $(H \o \iota_{M,N})\ci 
\bar{\iota}_{H, M\Co N}$, provided 
\begin{align}
(\mu_H \o M \o N)\,\ci \, \tw_{23} \ci (\bar{\iota}_{H, M} \o \bar{\iota}_{H, N})
\ci \iota_{H\Co M, H\Co N} ((h \o m) \o (h' \o n)) = \nonumber \\
= (h h') \o (m \o n) \in  H \coCe (M \coC N) \nonumber 
\end{align} 
This leads to the following equations:
\medskip
\begin{align}
(h h') \o m_{[0]} \o m_{[1]} \o n &=  
(h h') \o m \o n_{[-1]} \o n_{[0]} \label{egy}\\
(h h')_{(1)} \o \alpha ((h h')_{(2)}) \o m \o n &= (h h') \o m_{[-1]} \o 
m_{[0]} \o n \label{ketto}\\
(h h')_{(1)} \o \beta ((h h')_{(2)}) \o m \o n &= (h h') \o n_{[1]} \o m \o 
n_{[0]} \label{harom}
\end{align}

\medskip
Observe that by the multiplicativity of the coproduct and because
$(H, \lambda_L, \rho_L)$ is a monoid in ${\bicomC}$, we have the following identities:
\begin{align}
\alpha (h h') &= \alpha (h) \varepsilon (h') \label{elso}\\
\beta (h h') &= \varepsilon (h) \beta (h')\label{masodik}
\end{align}
To show \ref{elso}, compute
\begin{align}\nonumber
\alpha (h h')  = \alpha ((h h')_{(1)}) \varepsilon ((hh')_{(2)}) =
\alpha (h_{(1)}) \varepsilon (h_{(2)}) \varepsilon (h') = \alpha (h)
\varepsilon (h'),
\end{align}
and analagously for \ref{masodik}. Note that \ref{elso} and \ref{masodik} are dual to the relations
$\Delta_A (t(r)) = 1_A \o t(r)$ and $\Delta_A(s(r)) = s(r) \o 1_A$, which hold for a left
bialgebroid $A$ over $R$. 

\medskip
To prove \ref{ketto}, use \ref{elso} in the first equality and \ref{A}
in the second: 
\begin{align} 
(h h')_{(1)} \o \alpha ((h h')_{(2)}) \o m \o n &= h_{(1)} h' \o \alpha
(h_{(2)}) \o m \o n = \nonumber \\
&= (h h') \o m_{[-1]} \o m_{[0]} \o n \nonumber
\end{align}

Similarly, \ref{harom} is by proved by using \ref{masodik} in the
first equality and \ref{C} in the second: 
\begin{align} 
(h h')_{(1)} \o \beta ((h h')_{(2)}) \o m \o n &= h h'_{(1)} \o 
\beta (h'_{(2)}) \o m \o n = \nonumber \\
&= (h h') \o n_{[1]} \o m \o n_{[0]}, \nonumber
\end{align}

To prove \ref{egy}, applying \ref{B} and \ref{C} to the left-- and
right hand sides, respectively, yields
\beq\nonumber
h_{(1)} h' \o m \o \beta(h_{(2)}) \o n = h h'_{(1)} \o m \o \alpha
(h'_{(2)}) \o n
\eeq
which holds precisely because multiplication satisfies the property \ref{brzmil1}.

\medskip
Let us check that $\delta_{M,N}: M \coC N \to H \coCe (M \coC N)$ is indeed a coaction. Expressed
on elements, $\delta_{M,N} (m \o n) = m_{[ -1 ]} n_{[ -1 ]} \o m_{[ 0 ]}
\o n_{[ 0 ]}$ (we think of the domain and range of $\delta_{M,N}$ as embedded into $M \o N$ and 
$H \o (M \o N)$, respectively). 
\begin{align}
(H \Co \delta_{M,N})\ci \delta_{M,N} (m \o n) &= m_{[-1]} n_{[-1]} \o 
m_{[0] [-1]} n_{[0] [-1]} \o m_{[0]}
\o n_{[0]} = \nonumber\\
&= m_{[-1] (1)} n_{[ -1 ] (1)} \o  m_{[ -1 ] (2)} n_{[ -1 ] (2)}
\o m_{[ 0 ]}\o n_{[ 0 ]} = \nonumber\\
&= (m_{[ -1 ]} n_{[ -1 ]})_{(1)} \o
(m_{[ -1 ]} n_{[ -1 ]})_{(2)} \o m_{[ 0 ]}
\o n_{[ 0 ]} = \nonumber\\
&= (\Delta \o M) \ci \delta_{M,N} (m \o n) \nonumber
\end{align}
where we used the comultiplicativity of the multiplication on $H$ in
the third equality.

\medskip
For $\langle ^H\M,\Co, C \rangle$ to be a monoidal category, we have still to
define the natural isomorphisms $\alpha_{M,N,P}: (M \Co N) \Co P \to M
\Co (N \Co P)$ (the associator), $\lambda_M: C \Co M \to M$ and
$\rho_N: N \Co C \to N$. Due to the strict monoidality of $F$, these maps may be defined
as the lifting of the respective coherence morphisms of ${\bicomC}$ to $^H\M^H$, 
provided they induce $H$--comodule maps. This, however, follows from the associativity and unit property 
of the multiplication and unit on $H$.   
\end{proof}

\subsection{Modules over a bicoalgebroid}
We proceed to define modules over a bicoalgebroid, especially for the purposes of
Section 3. 

\begin{defi}
A right module over a left bicoalgebroid $H$ (over $C$) is a pair $\langle X, \ract \rangle$, where
$X \in {\M^C}$ is a right $C$--comodule and the action is a right $C$--comodule map 
$\ract: X \coC H^{C} \to X^C$.
Similarly, a left module is a pair $\langle Y, \lact \rangle$ with $Y \in {^C\M}$ and
$\lact: {^{C}H} \coC Y \to Y$ a left $C$--comodule map. $H$ is a
$C$--bicomodule through the coactions $\lambda_L$ and $\rho_L$.
\end{defi}

The module category of a bicoalgebroid is expected to be monoidal as well, coming with
an embedding into $\bicomC$. The above definition doesn't seem to allow for this, but luckily,
a dual of Prop.~1.1. of \cite{Kor-Bim} holds:

\begin{pro}
Let $\langle X, \ract \rangle$ be a right module over the bicoalgebroid $H$. Then $X$ has a
unique left $C$--comodule structure such that
\begin{enumerate}
\item $X$ is a $C$--bicomodule
\item the action is a $C$--bicomodule map
\item $\ract: X \coC H \to X$ factorizes through $X \boxtimes H$
\end{enumerate}
\end{pro}
\begin{proof}
Note that the action being a right $C$--comodule map means
\beq\label{rightcoac} (x \ract h)_{[0]} \o (x \ract h)_{[1]} = x\ract h\twoB \o \beta (h\oneB) 
\eeq 
The left comodule structure in question will be denoted $\tau (x) = x_{[-1]} \o x_{[0]}$. In fact,
$\tau$ is uniquely determined by demanding that the right $H$ action
be also a left $C$--comodule map w.r.t $\tau$. Note that $X \coC H$ is a left $C$--comodule
through the left $C$--coaction $\lambda_R (h) = \beta
(h_{(2)}) \o h_{(1)}$ on $H$, i.e.~we impose:
\beq\label{leftcoac} (x \ract h)_{[ -1 ]} \o (x \ract h)_{[ 0 ]} =
\beta (h_{(2)}) \o x\ract h_{(1)} 
\eeq 
The identity $x = x_{[ 0 ]} \ract \eta (x_{[ 1 ]})$ and \ref{leftcoac} yield an explicit formula
for the left coaction $\tau$:
\begin{align}
x_{[ -1 ]} \o x_{[ 0 ]} &= (x_{[ 0 ]} \ract \eta (x_{[ 1 ]}))_ {[ -1
] } \o (x_{[ 0 ]} \ract \eta (x_{[ 1 ]}))_
{[ 0 ] } = \nonumber \\
&=\beta (\eta (x_{[ 1 ]})_{(2)}) \o x_{[ 0 ]} \ract \eta
(x_{[ 1 ]})_{(1)}\nonumber
\end{align}
This is indeed a coaction, i.e.~$(C \o \tau) \ci \tau = (\Delta_C \o X) \ci \tau$. Inserting definitions, the
\begin{align}
LHS = \beta (\eta(x_{[ 1 ]})_{(2)}) &\o \beta\{\eta [(\xp)_{[ 1 ]}]_{(2)}\}
\o \nonumber \nonumber \\
& \o (\xp)_{[ 0 ]} \ract \eta [(\xp)_{[ 1 ]}]_{(1)} \nonumber
\end{align}
Using \ref{rightcoac}, we find:
\begin{align}
LHS = \beta (\eta(x_{[ 1 ]})_{(2)}) \o
\beta \{ \eta [\alpha(& \eta(x_{[ 1 ]})_{(1)(2)})]_{(2)}\} \o \nonumber \\
&\o (x_{[ 0 ]} \ract \eta(x_{[ 0 ]})_{(1)(1)}) \ract \eta
[\alpha (\eta(x_{[ 1 ]})_{(1)(2)})]_{(1)} \nonumber
\end{align}
which, by the bicoalgebroid axiom \ref{etaalphaeta}, is further equal:
\begin{align}
LHS &= \beta (\eta(x_{[ 1 ]})_{(2)}) \o \beta (\eta(x_{[ 1 ]})_{(1)(2)})_{(2)} \o
x_{[ 0 ]} \ract \eta(x_{[ 1 ]})_{(1)(1)} \eta(x_{[ 1 ]})_{(1)(2)(1)} = \nonumber\\
&= \beta (\eta(x_{[ 1 ]})_{(3)}) \o \beta (\eta(x_{[ 1 ]})_{(2)}) \o
x_{[ 0 ]} \ract \eta(x_{[ 1 ]})_{(1)} = \beta (\eta(x_{[ 1 ]})_{(2)})_{(1)} \o
\nonumber\\
& \o \beta (\eta(x_{[ 1 ]})_{(2)})_{(2)} \o x_{[ 0 ]} \ract \eta(x_{[ 1 ]})_{(1)}
= RHS .  \nonumber
\end{align}
In the first equality, we used comultiplicativity of the unit and coassociativity. In the second, the fact
that $\beta$ is an anti--coalgebra map. \medskip\\

As for (1), the coaction $\tau$ makes $X$ a bicomodule. Using the definition of the left coaction, and that
the $H$--action is a right $C$ comodule map:
\begin{align}
x_{[ -1 ]} &\o x_{[ 0 ]
[ 0 ]} \o x_{[ 0 ] [ 1 ]} =
\beta (\eta (x_{[ 1 ]})_{(2)}) \o  (x_{[ 0 ]} \ract
\eta(x_{[ 1 ]})_{(1)})_{[ 0 ]} \o \nonumber\\
&\o (x_{[ 0 ]} \ract
\eta(x_{[ 1 ]})_{(1)})_{[ 1 ]} = \beta (\eta (x_{[ 1 ]})_{(2)}) \o
x_{[ 0 ]} \ract \eta (x_{[ 1 ]})_{(1)(2)} \o \beta (\eta (x_{[ 1 ]})_{(1)(1)})\nonumber
\end{align}
Using that $\eta : {C^C} \to {H^C}$ is a $C$--bicomodule map,
\beq \label{etabikommap}  
\eta(c_{(1)}) \o c_{(2)} = \eta(c)_{(2)} \o \beta (\eta(c)_{(1)})
\eeq

and the coassociativity of the coaction:
\begin{align}\nonumber
x_{[ -1 ]} &\o x_{[ 0 ][ 0 ]} \o x_{[ 0 ] [ 1 ]} = 
\beta (\eta (x_{[ 1 ]})_{(2)(2)}) \o x_{[ 0 ]} \ract \eta (x_{[ 1 ]})_{(2)(1)} \o 
\beta (\eta(x_{[1]})_{(1)}) = \nonumber \\
&= \beta (\eta (x_{[1](1)})_{(2)}) \o x_{[0]} \ract \eta (x_{[1](1)})_{(1)} \o x_{[1](2)} = 
 \beta (\eta (x_{[0][1]})_{(2)}) \o \nonumber \\
&\o x_{[0][0]} \ract \eta (x_{[0][1]})_{(1)} \o x_{[1]} = x_{[0][-1]} \o  x_{[0][0]} \o x_{[1]} \nonumber
\end{align}
(we apply \ref{etabikommap} to $c = x_{[1]}$ in the second equality). 
The action will then (by construction) be a $C$--bicomodule map, proving (2). It remains to see that the action
factorizes through the cocenter of $X \coC H$, meaning:
\beq
(x_{[ 0 ]} \ract \eta(x_{[ 1 ]})_{(1)})\ract h \o \beta(\eta(x_{[ 1 ]})_{(2)}) =
x_{[ 0 ]} \ract h \o x_{[ -1 ]}
\eeq
This is a simple consequence of \ref{leftcoac}:
\begin{align}\nonumber
LHS = (x_{[ 0 ]} \ract \eta(x_{[ 1 ]}))_{[ 0 ]} \ract h \o
(x_{[ 0 ]} \ract \eta(x_{[ 1 ]}))_{[ -1 ]} = RHS.
\end{align}
\end{proof}

\section{The scalar extension for bicoalgebroids}
In \cite{Brz-Mil}, the authors introduced a construction that
associates to a bialgebra $H$ and a \emph{braided commutative
algebra} $Q$ over $H$ a bialgebroid. In \cite{Kor-Bim}, it was shown
that the construction generalizes to bialgebroids (in fact, even to
Frobenius Hopf--algebroids) and has an interpretation as the
noncommutative scalar extension of $H$ by $Q$. 

Here we dualize this construction to bicoalgebroids, and give a few simple examples. We begin by
defining the smash coproduct (\cite{SM:book}, with a slight variation).

\begin{defi}
Let $H$ be a bicoalgebroid over $C$ and $D$ an $H$--comodule coalgebra. Then
their smash coproduct $D \shrp H$ is a coalgebra, isomorphic to $D \coC H$ as
$C$--bicomodules and with the coalgebra structure:
\begin{align}
\Delta (d \shrp h) &= d_{(1)} \shrp {d_{(2)}}\coa{-1} h_{(1)} \coD {d_{(2)}}\coa{0}
\shrp h_{(2)}\\
\varepsilon (d \shrp h) &= \varepsilon (d) \varepsilon_H (h)
\end{align}
\end{defi}
That these maps define a coalgebra is easily verified. The category
of $(D \shrp H)$--comodules may also be described as the internal
$D$--comodules in $^H\M$, i.e.~$^D{(^H\M)} = {^{D \shrp H}\M}$.

Indeed, assume $X \in {^D{(^H\M)}}$. To every coaction $\delta_D
: X \to D \coC X$ in $^H\M$, we can associate a coaction of $D \shrp H$,
namely $\delta_{D\shrp H} = (D \o \delta) \ci \delta_D: X \to D \o X
\to D \o (H \o X)$, $\delta_{D\shrp H} (x) = x_{[-1]} \o {x_{[0]}}\coa{-1} \o {x_{[0]}}\coa{0}$.
A straightforward calculation proves that $(\Delta_{D\shrp H} \o X) \ci \delta_{D\shrp H} =
((D\shrp H) \o \delta_{D\shrp H}) \ci \delta_{D\shrp H}$, using that $\delta_{D\shrp H}$ is an $H$--comodule map.
In the reverse direction, an $(D \shrp H)$--comodule is both an $H$-comodule and a $D$--comodule such that the
$D$--coaction is an $H$--comodule map, which means precisely that it is an internal $D$--comodule in $^H\M$.

\subsection{Cocommutative coalgebras over bicoalgebroids}
Keeping with the method of reversing arrows, we arrive at the
following definition for Yetter--Drinfel'd modules over a
bicoalgebroid.

\begin{defi}
Let $H$ be a (left--) bicoalgebroid over $C$. A Yetter--Drinfel'd
module over $H$ is a triple $\langle Z, \ract, \delta \rangle$ such
that the $C$--bicomodule $Z$ is simultaneously a right $H$--module
with $\ract: Z \coC H \to Z$ and a left $H$--comodule with $\delta: Z
\to H \coCe Z$ so that the action and coaction satisfy the
compatibility condition

\beq\label{YD}
d\coa{-1} h\oneL \coC d\coa{0} \ract h\twoL = h\twoL
    (d \ract h\oneL)\coa{-1} \coC (d \ract h\oneL)\coa{0}
\eeq
\end{defi}

The Yetter--Drinfel'd category, denoted ${^H\YD _H}$ over $H$ has objects the
Yetter--Drinfel'd modules over $H$ and arrows the $C$--bicomodule maps
that are at the same time $H$--module maps and $H$--comodule maps.

\medskip
The category ${^H\YD _H}$ becomes monoidal if we define the monoidal
product of two Yetter--Drinfel'd modules $Z$, $Z'$ as $Z \coC Z'$ with
action and coaction:
\begin{align}
(z \coC z') \ract h &= z \ract h_{(2)} \coC z' \ract h_{(1)} \nonumber\\
(z \coC z')_{\langle -1 \rangle} \coC (z \coC z')_{\langle 0 \rangle}
&= z_{\langle -1 \rangle} z'_{\langle -1 \rangle} \coC  z_{\langle 0
  \rangle} \coC z'_{\langle 0 \rangle} \nonumber
\end{align}
The monoidal unit is of course $C$, with $c \ract h = c\, \varepsilon
(h)$ and $c_{\langle -1 \rangle} \o c_{\langle 0 \rangle} = \eta
(c_{(1)}) \o c_{(2)}$. Moreover, ${^H\YD _H}$ is pre--braided with
\beq
\tau_{Z,Z'}: Z \coC Z' \to Z' \coC Z, \;\; z \o z' \mapsto
z'_{\langle 0 \rangle} \coC z \ract z'_{\langle -1 \rangle}
\eeq

\bigskip

From experience with Hopf algebras, weak Hopf algebras and bialgebroids, it is reasonable to expect that 
the Yetter--Drinfel'd category over a bicoalgebroid is related to the (weak) center of
the category of comodules. For the center construction, consult \cite{Majid}, \cite{Joyal-Street2} and 
\cite{Kassel}. The notion of weak center of a monoidal category seems to appear in    
\cite{Schauenburg: ddqg}, Definition 4.3 (see also \cite{Cae-Wang-Yin}, Section 1.3 and \cite{Kor-Bim}).

For bialgebroids, it is known that the Yetter--Drinfel'd category is equivalent to the monoidal weak 
center (see \cite{Schauenburg: ddqg}). Unfortunately this doesn't seem to be true for bicoalgebroids in 
general. Nevertheless, the $\YD$ category over a bicoalgebroid still embeds into the monoidal weak center. 
Although the weak center construction is applicable to any monoidal category, we shall only recall the 
definition in the context of the comodule category over a bicoalgebroid.

For a bicoalgebroid $H$ over $C$, the (left) weak center $\lZ({^H\M})$ has
objects $\langle Z, \theta \rangle$, where $Z \in {^H\M}$ and $\theta$ is a natural
transformation $\theta_Y : Z \coC Y \to Y \coC Z$ (between endofunctors on ${^H\M}$)
that satisfies
\begin{align}
\theta_{X \coC Y} &= (X \coC \theta_Y) \ci (\theta_X \coC Y) \label{centrum1}\\
\theta_C &= Z \label{centrum2}
\end{align}
An arrow $\langle Z, \theta \rangle \to \langle Z', \theta' \rangle$ is an $H$--comodule map
$f: Z \to Z'$, compatible with $\theta$'s in the sense:
\beq
(Y \coC f) \ci \theta_Y = \theta'_Y \ci (f \coC Y)
\eeq
for all $Y \in {^H\M}$. The category $\lZ({^H\M})$ is monoidal and pre--braided with
monoidal product
\beq
\langle Z, \theta \rangle \coC \langle Z', \theta' \rangle = \langle Z
\coC Z', (\theta\_\, \coC Z')\ci
(Z \coC \theta'\_\,) \rangle
\eeq
and pre--braiding
\beq
\lbeta_{\langle Z, \theta \rangle, \langle Z', \theta' \rangle} = \theta_{Z'}
\eeq

\bigskip
It is easily shown that every Yetter--Drinfel'd module $\langle Z,
\delta, \ract \rangle$ has the structure of an object in
$\lZ(^H\M)$. The map

\begin{align}
\theta_X : Z \coC X &\to X \coC Z \\
z \o x &\mapsto x\coa{0} \o z \ract x\coa{-1} \nonumber
\end{align}
is natural in $X$, since the arrows of $\lZ(^H\M)$ are $H$--comodule
maps, $\theta_C = Z$ is trivially satisfied and
\beq\nonumber
\theta_{X \Co Y} (x
\o y) = x\coa{0} \o y\coa{0} \o z\ract (x\coa{-1} y\coa{-1})
\eeq
equals
\begin{align}
(X \Co \theta_Y)\ci (\theta_X \Co Y) (x \o y) &= (X \Co \theta_Y)
(x\coa{0} \o z\ract x\coa{-1} \o y) = \nonumber\\
&= x\coa{0} \o y\coa{0} \o (z\ract
x\coa{-1}) \ract y\coa{-1}. \nonumber
\end{align}

\bigskip

As for the reverse direction, we can associate to every object $\langle Z, \theta \rangle$ of $\lZ({^H\M})$
a right action of $H$ as follows:
\begin{align}
\ract&: Z \coC H \to Z\\
&z \o h \mapsto (\varepsilon_H \o Z) \ci \theta_H (z \o h)\nonumber
\end{align}
It is easily checked that this is indeed a right $H$--action, and is 
the candidate to make $Z$ a Yetter--Drinfel'd module. If $\theta$ enjoys the property
$\theta_X (z \o x) = x\coa{0} \o \theta_H(z \o x\coa{-1})$ for all
objects $X \in {^H\M}$, then $\langle Z, \delta, \ract \rangle$ becomes a
Yetter--Drinfel'd module and, moreover, $^H\YD_H$ and $\lZ(^H\M)$ are
isomorphic categories. This would mean that the natural map $\theta$
can be expressed with it's component $\theta_H$. This is indeed
possible for \emph{bialgebroids}, since any bialgebroid is a generator
in the category of modules over itself, and natural transformations are
determined by their value on the generator.

\medskip
\begin{rmk}
We mention, for the sake of completeness, the \emph{right} weak center $\rZ(^H\M)$, defined as the category of
pairs $\langle Z, \bar{\theta} \rangle$, where $Z \in {^H\M}$ and $\bar{\theta}$ is a natural
transformation $\theta_Y : Y \coC Z \to Z \coC Y$, satisfying
\begin{align}
\bar{\theta}_{X \coC Y} &= (\bar{\theta}_X \coC Y) \ci (X \coC \bar{\theta}_Y) \label{rcentrum1}\\
\theta_C &= Z \label{rcentrum2}
\end{align}
The category $\rZ({^H\M})$ has the monoidal structure
\beq
\langle Z, \bar{\theta} \rangle \coC \langle Z', \bar{\theta}' \rangle = \langle Z
\coC Z', (Z \coC \bar{\theta}'\_\,) \ci
(\bar{\theta}\_\,\coC Z') \rangle
\eeq
and pre--braiding
\beq
\rbeta_{\langle Z, \bar{\theta} \rangle, \langle Z', \bar{\theta}' \rangle} = \bar{\theta}'_{Z}
\eeq

\medskip
It is straightforward to prove that the one--sided Yetter--Drinfel'd category $^H_H\YD$ is embedded into
$\rZ(^H\M)$. The objects of $^H_H\YD$ are triples $\langle Z, \delta, \lact \rangle$, $C$--bicomodules which are
simultaneously $H$--modules and $H$--comodules, satifying the compatibilty condition
\beq
h_{(1)} z\coa{-1} \coC h_{(2)} \lact z\coa{0} = (h_{(1)} \lact z)\coa{-1} h_{(2)} \coC (h_{(1)} \lact z)\coa{0}
\eeq
$^H_H\YD$ is a pre--braided monoidal category with the pre--braiding
\beq
\kappa_{Z',Z} (z' \o z) = {z'}\coa{-1} \lact z \o {z'}\coa{0}
\eeq
\end{rmk}

\bigskip
Now, a braided cocommutative coalgebra (hereinafter abbreviated BCC)
over $H$ is defined as a cocommutative comonoid in
${^H\YD _H}$. Spelled out in detail, we have the

\begin{defi}
A BCC over $H$ is a coalgebra $D$, equipped with a coalgebra map
$\varepsilon: D \to C$ and the structure of a Yetter--Drinfel'd module
$\langle D, \ract, \delta \rangle \in {^H\YD _H}$ so that the
left/right $C$--comodule structures on $D$ are given by $\varepsilon
(d_{(1)}) \o d_{(2)}$ and $d_{(1)} \o \varepsilon(d_{(2)})$,
respectively and the relations stating that $D$ is an $H$--module and $H$--comodule
coalgebra:
\begin{align}
(d \ract h)_{(1)} \o (d \ract h)_{(2)} &= d_{(1)} \ract h_{(1)} \o
d_{(2)} \ract h_{(2)}\\
\varepsilon (d \ract h) &= \varepsilon (d) \varepsilon_H (h)\\
{d_{(1)}}\coa{-1} {d_{(2)}}\coa{-1} \o {d_{(1)}}\coa{0} \o {d_{(2)}}\coa{0} &= d\coa{-1} \o
  {d\coa{0}}_{(1)} \o {d\coa{0}}_{(2)}\\
d\coa{-1} \o \varepsilon (d\coa{0}) &=
\eta (\varepsilon (d)_{(1)}) \o \varepsilon (d)_{(2)}\label{BCC4}
\end{align}
and braided cocommutativity:
\beq\label{brcocom}
d_{(1)} \o d_{(2)} = {d_{(2)}}\coa{0} \o d_{(1)}\ract {d_{(2)}}\coa{-1}
\eeq
\end{defi}

We have the following functorial characterization of BCC's, entirely
analogous to Prop 4.7. of \cite{Kor-Bim}:

\begin{lem}\label{BCClem}
If $D$ is a BCC in ${^H\YD _H}$, then the functor $D \coC \_: {^H\M}
\to {^D{(^H\M)}} = {^{D \shrp H}\M}$ is strong monoidal
\end{lem}

\begin{proof}
Denote the opmonoidal structure $\langle {D \coC \_\,}, D^2, D^0
\rangle$. The natural transformation
\begin{align}
(D^2)_{X,Y} : &D \coC (X \coC Y) \to (D \coC X) \coD (D \coC Y)\nonumber\\
& d \o x \o y \mapsto (d_{(1)} \o x_{\langle 0 \rangle}) \o
(d_{(2)}\ract x_{\langle -1 \rangle} \o y)\nonumber
\end{align}
has the inverse $(d \o x) \o (d' \o y) \mapsto d \varepsilon_C(\varepsilon (d')) \o x
\o y$. Furthermore, $D^0: D \coC C \to D$ is obviously an
isomorphism.
\end{proof}

It is perhaps not altogether surprising that we have the following
dualization of Theorem 4.6. of \cite{Kor-Bim}

\begin{thm}
Let $\langle H, \Delta, \varepsilon ; \mu, \eta ; \alpha , \beta ; C
\rangle$ be a (left--) bicoalgebroid over $C$ and $D$ a BCC over $H$,
then $\langle D \shrp H, \tilde{\Delta}, \tilde{\varepsilon} ; \tilde{\mu},
\tilde{\eta} ; \tilde{\alpha} , \tilde{\beta} ; D \rangle $ is a (left--) bicoalgebroid over $D$, with
the following structure maps:

\begin{align}
\tilde{\Delta} (d \shrp h) &= d_{(1)} \shrp {d_{(2)}}\coa{-1} h_{(1)} \coD {d_{(2)}}\coa{0}
\shrp h_{(2)}\\
\tilde{\varepsilon} (d \shrp h) &= \varepsilon_C(\varepsilon (d)) \varepsilon_H (h) \\
\tilde{\mu} (d \shrp h \coD d' \shrp h') &= d \varepsilon_C(\varepsilon (d')) \shrp h h' \\
\tilde{\eta} (d) &= d_{(1)} \shrp \eta(\varepsilon(d_{(2)})) \\
\tilde{\alpha} (d \shrp h) &= d \varepsilon_H (h), \;\;\; \tilde{\beta} (d \shrp h) = d \ract h
\end{align}
\end{thm}

\begin{proof}
First, we check that $\tilde{\alpha}$ ($\tilde{\beta}$) is a coalgebra
(anti--coalgebra) map, respectively. Inserting the definitions, a trivial calculation shows
\beq\nonumber
\alpha ((d\shrp h)_{(1)}) \o \alpha ((d\shrp h)_{(2)}) = d_{(1)} \o d_{(2)} \varepsilon (h) =
(\alpha(d \shrp h))_{(1)} \o (\alpha(d \shrp h))_{(2)}
\eeq
As required, $\beta$ is an anti--coalgebra map:
\begin{align}
&\beta ((d\shrp h)_{(2)}) \o \beta ((d\shrp h)_{(1)}) = \beta ({d_{(2)}}\coa{0} \shrp h_{(2)}) \o
\beta (d_{(1)} \shrp {d_{(2)}}\coa{-1} h_{(1)}) = \nonumber\\
&= {d_{(2)}}\coa{0} \ract h_{(2)} \o d_{(1)} \ract {d_{(2)}}\coa{-1} h_{(1)} = d_{(1)}\ract h_{(2)} \o
d_{(2)}\ract h_{(1)} = (d\ract h)_{(1)}\o (d\ract h)_{(2)} = \nonumber\\
&= (\beta (d\shrp h))_{(1)} \o (\beta (d\shrp h))_{(2)}. \nonumber
\end{align}
where we have used \ref{brcocom} in the third equality, and the fact that $D$ is an $H_{cop}$--coalgebra
in the fourth.

To prove that $\tilde{\mu}: (D \shrp H) \coD (D \shrp H) \to D \shrp H$ factorizes through
$(D \shrp H) \boxtimes (D \shrp H)$, we calculate the $D$--comodule structure of $D \shrp H$:
\begin{align}
\tilde{\lambda}_L &: d\shrp h \mapsto \alpha ((d\shrp h)_{(1)}) \o (d\shrp h)_{(2)} =
\alpha(d_{(1)} \shrp {d_{(2)}}\coa{-1} h_{(1)}) \o {d_{(2)}}\coa{0}
\shrp h_{(2)}= \nonumber\\
&= d_{(1)} \o d_{(2)} \shrp h \nonumber
\end{align}
\begin{align}
\tilde{\rho}_L &: d\shrp h \mapsto (d\shrp h)_{(2)} \o \beta((d\shrp h)_{(2)}) =
{d_{(2)}}\coa{0} \shrp h_{(2)} \o d_{(1)} \ract ({d_{(2)}}\coa{-1} h_{(1)}) =\nonumber \\
&= d_{(1)} \shrp h_{(2)} \o d_{(2)} \ract h_{(1)}, \nonumber
\end{align}
using (\ref{brcocom}) in the last step. The definition of the cotensor product over $D$ then reads:
$(d\o h) \o (d'\o h') \in (D\shrp H) \coD (D\shrp H)$ \emph{iff}
$$(d\shrp h)_{(2)} \o \beta ((d\shrp h)_{(1)}) \o (d' \shrp h') = (d\shrp h) \o \alpha
((d'\shrp h')_{(1)}) \o (d' \shrp h')_{(2)}$$
or, using (\ref{brcocom}):
\beq\label{cotD}
d_{(1)} \shrp h_{(2)} \o d_{(2)} \ract h_{(1)} \o d' \shrp h' =
d\shrp h \o d'_{(1)} \o d'_{(2)} \shrp h'
\eeq

We now prove $(d\shrp h)(d'\shrp h')_{(1)} \o \alpha((d'\shrp h')_{(2)}) = (d\shrp h)_{(1)}(d'\shrp h') \o
\beta ((d\shrp h)_{(2)})$. Inserting definitions, and using the Yetter--Drinfel'd condition (\ref{YD}) we find:
\begin{align}
RHS &= d_{(1)}\varepsilon_C(\varepsilon(d')) \shrp {d_{(2)}}\coa{-1}h_{(1)}h' \o {d_{(2)}}\coa{0} \ract h_{(2)} = 
\nonumber\\
 &= d_{(1)}\varepsilon_C(\varepsilon(d')) \shrp h_{(2)} (d_{(2)} \ract h_{(1)})\coa{-1} h' \o 
(d_{(2)}\ract h_{(1)})\coa{0},
\nonumber
\end{align}
using the Yetter--Drinfel'd condition (eq.~\ref{YD}). Applying (\ref{cotD}), we arrive at
\begin{align}
RHS = d \varepsilon_C(\varepsilon (d'_{(2)})) \shrp h {{d'}_{(1)}}\coa{-1}h' \o {{d'}_{(1)}}\coa{0} = d \shrp h
{d'}\coa{-1} h' \o {d'}\coa{0}\nonumber
\end{align}
A quick calculation shows that the
\begin{align}
LHS = d \varepsilon_C(\varepsilon ({d'}_{(1)})) \shrp h {{d'}_{(2)}}\coa{-1} h' \o {{d'}_{(2)}}\coa{0} = d \shrp
h {d'}\coa{-1} h' \o {d'}\coa{0},\nonumber
\end{align}
as claimed.

Comultiplicativity of the product (which makes sense due to our above assertion) means
\beq\nonumber
(\tilde{\Delta} \ci \tilde{\mu}) [(d\shrp h) \coD (d'\shrp h')] =
(\tilde{\mu} \coD \tilde{\mu})\ci \tau_{23} \ci (\tilde{\Delta} \coD \tilde{\Delta}) [(d\shrp h) \coD (d'\shrp h')]  
\eeq
inserting our definitions, we have:
\beq\nonumber
LHS = \tilde{\Delta} (d\varepsilon_C\varepsilon (d') \shrp hh') = 
d_{(1)}\varepsilon_C\varepsilon (d') \shrp {d_{(2)}}\coa{-1} (hh')_{(1)} \coD {d_{(2)}}\coa{0}\shrp (hh')_{(2)},
\eeq
on the other hand, the
\begin{align}\nonumber
RHS &= (d_{(1)}\shrp {d_{(2)}}\coa{-1}h_{(1)}) (d'_{(1)}\shrp {d'_{(2)}}\coa{-1}h'_{(1)}) \coD 
({d_{(2)}}\coa{0}\shrp h_{(2)}) ({d'_{(2)}}\coa{0}\shrp h'_{(2)})=\\
&= d_{(1)}\varepsilon_C\varepsilon(d'_{(1)})\shrp {d_{(2)}}\coa{-1}h_{(1)} {d'_{(2)}}\coa{-1}h'_{(1)} \coD
{d_{(2)}}\coa{0}\varepsilon_C\varepsilon({d'_{(2)}}\coa{0}) \shrp h_{(2)}h'_{(2)} =\nonumber\\
&= d_{(1)}\varepsilon_C(\varepsilon(d'_{(1)}))\shrp {d_{(2)}}\coa{-1}h_{(1)} \eta(\varepsilon(d'_{(2)}))h'_{(1)} \coD
{d_{(2)}}\coa{0} \varepsilon_C(\varepsilon(d'_{(3)})) \shrp (hh')_{(2)}, \nonumber
\end{align}
where we made use of \ref{BCC4} and coassociativity in the third
equality. Now, $d'_{(1)} \o \varepsilon(d'_{(2)}) \o
h' = d' \o \alpha(h'_{(1)}) \o h'_{(2)}$, because $d\shrp h \in D\coC
H$. From this, and the unit property of $\eta$, the statement follows.

The product is counital:
\begin{align}
\tilde{\varepsilon} (d\shrp h) \tilde{\varepsilon}(d'\shrp h') =
\varepsilon_C(\varepsilon (d))\varepsilon_C(\varepsilon (d'))
\varepsilon_H (hh') = \tilde{\varepsilon} (d\varepsilon_C(\varepsilon
(d)) \shrp hh')  
\end{align}

The unit map $\tilde{\eta}$ is indeed a unit for $\tilde{\mu}$. The first unit 
property reads:
\begin{align}
&\tilde{\mu}\ci (\tilde{\eta} \Co D\shrp H)\ci \tilde{\lambda}_L
(d\shrp h) = (d_{(1)(1)} \shrp \eta(\varepsilon (d_{(1)(2)})))
(d_{(2)}\shrp h) = \nonumber\\
&=d_{(1)}\varepsilon_C\varepsilon(d_{(3)}) \shrp
\eta(\varepsilon(d_{(2)})) h = d\shrp h,
\end{align}
using $d\shrp h \in D\coC H$ and the unit axiom (for $H$) in the last
equality. The second,
\begin{align}
&\tilde{\mu}\ci (D\shrp H \Co \tilde{\eta})\ci \tilde{\rho}_L
(d\shrp h) = (d_{(1)}\shrp h_{(2)})  ((d_{(2)}\ract h_{(1)})_{(1)}
\shrp \eta((d_{(2)}\ract h_{(1)})_{(2)})) = \nonumber\\
&= d_{(1)} \varepsilon_C\varepsilon (d_{(2)}\ract h_{(1)}) \shrp
h_{(2)} \eta (\varepsilon(d_{(2)} \ract h_{(1)})_{(2)}) = d_{(1)}
\shrp \eta(\varepsilon(d_{(2)})) h = d\shrp h 
\nonumber 
\end{align} 
is proved using that $D$ is an $H_{cop}$--algebra in the third
equality, and $d\shrp h \in D\coC H$ in the last. As a coalgebra, $D\shrp H$ is the
smash coproduct. The algebra structure of $\langle D\shrp H, \tilde{\mu}, \tilde{\eta}\rangle$ 
and the remaining axioms are easily verified.    
\end{proof}

\begin{exa}{\it{The action groupoid}}

\medskip
In the category $\Set$, there is a unique comultiplication, namely
the diagonal coproduct: $x \in X$, $\Delta_X(x) = x \times x$. The
counit is just a constant map to a (the) one--element set $1$, hence
$\varepsilon_X (x) = *$ for all $x \in X$, where $*$ is the unique element 
of $1$. The coaction of a group $G$ on $X$ is completely specified by an
arbitrary function $\varphi : X \to G$, via $\delta_{\varphi}(x) =
x_{\langle -1 \rangle} \times x_{\langle 0 \rangle} = \varphi (x)
\times x$. Now, consider a $G$--Set $\langle X, \ract \rangle$,
carrying a right action of $G$. Choosing a $G$--coaction
$\delta_{\varphi}$, the Yetter--Drinfel'd compatibility condition
takes the form \beq\label{YD-G} \varphi(x) g \times x\ract g = g
\varphi(x\ract g) \times x\ract g \eeq so $\langle X,
\delta_{\varphi}, \ract \rangle$ is a $\YD$--module in $^G\YD_G$ if
and only if $g^{-1} \varphi(x) g = \varphi (x\ract g)$. Moreover,
$X$ is a BCC if $x \times x = x \times x\ract \varphi (x)$, i.e.~iff
\beq\label{BCC-G} x \ract \varphi(x) = x \eeq \ref{BCC-G} implies
that the value of $\varphi$ at a point $x$ must lie in the
stabilizer subgroup $G^x$ of the point $x$, and from \ref{YD-G} we
conclude that it suffices to define $\varphi$ for a single
representative, say $x_0$ of each $G$-orbit. Then, if $x_0 \in
G^{x_0}$, $\varphi (x) = \varphi (x_0\ract g) = g^{-1} \varphi(x_0)
g \in G^x$.

\medskip
Choosing a trivial coaction $\varphi (x) \equiv e$, the scalar extension of $G$ by $X$ is nothing but the
action groupoid. Indeed, $\tilde{\alpha}(x \shrp g) = x$ and $\tilde{\beta}(x \shrp g) = x\ract g$, so
$(X \shrp G) \Co_X (X \shrp G)$ is the set of composable pairs in the action groupoid and the
multiplication $\tilde{\mu}$ is the composition of arrows in the action groupoid.

\medskip
The phenomenon behind this example is that in $\Set$, the fibered
product of two parallel maps $\alpha, \beta: X \to Y$, defined by
the pullback
$$\xymatrix{
X \times_{\alpha,\beta} X \ar[r]^-{q} \ar[d]_-{p}& X \ar[d]^-{\alpha}\\
X \ar[r]_{\beta} & Y }$$ is equivalent to the equalizer
$$\xymatrix{X \times_{\alpha,\beta} X \to X \times X \ar@<0.5ex>[r]^-{X \times \lambda}
\ar@<-0.5ex>[r]_-{\rho \times X} & X \times Y \times X}$$ where
$\lambda$ and $\rho$ are the 'coactions' $\lambda = (\alpha \times
X)\ci \Delta_{diag}$ and $\rho = (X \times \beta)\ci \Delta_{diag}$.
It is in this sense that a groupoid may be regarded as a classical
ancestor of a bicoalgebroid.
\end{exa}

\begin{exa}{\it{The regular BCC for $H$ a Hopf algebra}}

\medskip
$k$--Hopf algebras (and bialgebras) are examples both of bialgebroids and bicoalgebroids. It is not immaterial,
however whether we consider the Yetter-Drinfel'd category $^H\YD_H$ as embedded in $\lZ(\M_H)$ (the 'bialgebroid view',
see \cite{Kor-Bim}), or in $\lZ(^H\M)$ (the 'bicoalgebroid view'). Namely, the braiding is different in the two
cases, $\lZ(\M_H)$ is pre--braided with $\lbeta_{Z,Z'} = z'\ract z\coa{-1}
\o z\coa{0}$ and $\lZ(^H\M)$ is pre--braided with $\lgamma_{Z,Z'} =
{z'}\coa{0} \o z\ract {z'}\coa{-1}$.

\bigskip
A $k$--Hopf algebra $H$, with invertible antipode is a Yetter--Drinfel'd module $\langle H,
Ad_R, \Delta \rangle$ in  $\lZ(\M_H)$ (the \emph{regular} module) via
the coproduct, considered as left $H$--coaction and the right adjoint action,
$ Ad_R: H \o H \to H,\;\; h \o h' \mapsto S^{-1} (h'_{(2)})h
h'_{(1)}$. Furthermore, $H^{op}$ is a BCA in $\lZ(\M_H)$, that is $\mu_{op} \ci \beta (h \o
h') = \mu_{op} (h \o h')$. Indeed, $h'h = h'\ract h_{(1)} \o h_{(2)} = S(h_{(1)(2)})h'h_{(1)(1)} \o
h_{(2)}$.

\medskip
Dually, a $k$--Hopf algebra $H$ is a Yetter--Drinfel'd module $\langle H,
\mu, \tilde{Ad}_L \rangle$ in  $\lZ(^H\M)$ (the \emph{regular} module) via
the multiplication considered as a right action, and the left adjoint
coaction,
\begin{align}
\tilde{Ad}_L: H &\to H \o H \\
h &\mapsto S^{-1} (h_{(3)}) h_{(1)} \o h_{(2)} \nonumber
\end{align}
Yetter--Drinfel'd compatibility is easily checked:
\begin{align}
&h'_{(2)} (hh'_{(1)})\coa{-1} \o (hh'_{(1)})\coa{0} = h'_{(2)} S^{-1}((hh'_{(1)})_{(3)}) (hh'_{(1)}) \o
(hh'_{(2)}) = \nonumber \\
&=h'_{(2)} S^{-1}(h'_{(1)(3)}) S^{-1}(h_{(3)}) h'_{(1)(1)}h_{(1)} \o h_{(2)}h'_{(1)(2)}
= S^{-1}(h_{(3)})h_{(1)}h'_{(1)} \o h_{(2)}h'_{(2)} = \nonumber \\
&= h\coa{-1} h'_{(1)} \o h\coa{0}\ract h'_{(2)} \nonumber
\end{align}
As one might expect from the previous example, $H_{cop}$ is a BCC in
$\lZ(^H\M)$,
\beq\nonumber
\lbeta \ci \Delta_{cop} (h) = \lbeta (h_{(2)} \o h_{(1)}) = h_{(1)(2)} \o h_{(2)}S^{-1}(h_{(1)(3)})h_{(1)(1)}
= h_{(2)} \o h_{(1)}
\eeq

\medskip
To construct an example which does not require the invertibility of the antipode, consider
$^H_H\YD$ as being in the \emph{right} weak center $\rZ(^H\M)$. The Yetter--Drinfel'd condition takes the form
\beq
h_{(1)} z\coa{-1} \o h_{(2)} \lact z\coa{0} = (h_{(1)} \lact z)\coa{-1} h_{(2)} \o (h_{(1)} \lact z)\coa{0},
\eeq
and the pre--braiding is $\rbeta_{Z',Z}: {z'}\coa{-1} \lact z \o {z'}\coa{0}$. We find that
for an arbitrary Hopf algebra, $\langle H, Ad_L, \mu_H \rangle$ is a BCC in $^H_H\YD$, where
\begin{align}
Ad_L: H &\to H \o H \nonumber \\
h &\mapsto h_{(1)}S(h_{(3)}) \o h_{(2)}, \nonumber
\end{align}
and $\rbeta_{H,H}\ci \Delta (h)= h_{(1)(1)} S (h_{(1)(3)}) h_{(2)} \o h_{(1)(2)} = h_{(1)} \o h_{(2)}$.
\end{exa}

\section{The scalar extension as a comonad}
In this section, we give a (co--)monadic characterization of bicoalgebroids which can be seen as
dual to the results obtained for bialgebroids in \cite{Sz:EilMoor}. We also give a categorical
description of the bialgebroid and bi\emph{co}algebroid scalar extensions in terms of bimonads, and
bicomonads, respectively.

Recall that for a bicoalgebroid $H$, the forgetful functor $F: {^H\M} \to {^{C^e}\M}$ is strong
monoidal, and is \emph{left} adjoint to the induction functor $I = {^HH} \coCe \_\,: {^{C^e}\M}\to {^H\M}$.
By the standard Eilenberg--Moore construction (see \cite{MacLane}),
the adjunction $F \dashv I$ gives rise to a monad $\mathbb{T} = \langle T, \mu, \eta
\rangle$ on the category ${^H\M}$ with underlying endofunctor $T = I F: {^H\M}\to {^H\M}$ (
monad multiplication is $\mu = I \varepsilon F: T T \to T$, monad unit $\eta: {^H\M}
\to T$ is the unit of the adjunction) and a comonad $\mathbb{G} = \langle G, \Delta,
\varepsilon \rangle$ on the category ${^{C^e}\M}$ with underlying endofunctor $G = F I: {^{C^e}\M} \to
{^{C^e}\M}$ (comonad comultiplication is $\Delta = F \eta I: G \to G G$, counit $\varepsilon: G
\to {^{C^e}\M}$ is the counit of the adjunction). Denote $^{\mathbb{G}}\M$ the Eilenberg--Moore category 
of $\mathbb{G}$--coalgebras, then $^{\mathbb{G}}\M$ can be identified with ${^H\M}$, since 
$G = {^{C^e}H} \coCe \_\,$. Also, the canonical forgetful functor $F_G: {^{\mathbb{G}}\M} \to{^{C^e}\M}$ can
be identified with $F: {^H\M} \to {^{C^e}\M}$.  

\medskip
By Prop.~2.1. of \cite{Sz:EilMoor}, the (strong) opmonoidal structure on $F$ implies a monoidal
structure on the right adjoint $I$, and the adjunction is in the category of monoidal categories.
This implies that the unit and counit are monoidal natural transformations. The following definition is 
tailor--made (see \cite{Moer}):
\begin{defi}
Let $\langle \M, \Co, I \rangle$ be a monoidal category. Then a bicomonad on $\M$ is
a comonoid in the category of monoidal endofunctors from $\M$ to $\M$.
Thus, it is an endofunctor $G : \M \to \M$, furnished with:
\begin{itemize}
 \item a natural transformation $\kappa_{X,Y} : (G\, X) \Co (G\, Y) \to G\, (X \Co Y)$, and
 \item an arrow $\xi: C \to G\, C$
\end{itemize}
such that $\langle G, \kappa_{X,Y}, \xi \rangle$ is a monoidal
functor;
\begin{itemize}
 \item a natural transformation $\delta_X: G\, X \to G G X$ and
 \item a natural transformation $\varepsilon_X : G\, X \to X$
\end{itemize}
such that $\langle G, \delta, \varepsilon \rangle$ is a comonoid in
${\M}^{\M}$, and four compatibility axioms stating that $\delta$
is monoidal,
\begin{align}
\delta_{X \otimes Y} \ci \kappa_{X,Y} &= (G\,\kappa_{X,Y} \ci \kappa_{G\,X,G\,Y}) 
\ci (\delta_X \otimes \delta_Y) \label{a}\\
\delta_I \ci \xi &= G\,\xi \ci \xi \label{b}
\end{align}
and that $\varepsilon$ is monoidal
\begin{align}
\varepsilon_{X \otimes Y} \ci \kappa_{X,Y} &= \varepsilon_X \otimes
\varepsilon_Y \label{c}\\
\varepsilon \ci \xi = I \label{d}
\end{align}
\end{defi}

\begin{pro}
The endofunctor $G = F I = {^{C^e}H} \coCe \_\, : {^{C^e}\M} \to
{^{C^e}\M}$ is a monoidal comonad with the structure maps:
\begin{align}
\delta_X  &: H \coCe X \to H \coCe (H \coCe X) \\
          &h \o x \mapsto h_{(1)} \o (h_{(2)} \o x) \\
\varepsilon_X  &: H \coCe X \to X \\
               &h \o x \mapsto \varepsilon_H (h) x \\
\kappa_{X,Y} &: (H \coCe X) \coC (H \coCe Y) \to H \coCe (X \coC Y) \\
             &(h \o x) \o (h' \o y) \mapsto hh' \o (x \o y)\\
\xi  &: C \to H \coCe C \\
       &c \mapsto \eta(c)_{(1)} \o \alpha (\eta(c)_{(2)})
\end{align}
\end{pro}

\begin{proof}
The associativity of $\kappa$ corresponds to the associativity of the
multiplication $\mu$ of $H$, and $\xi$ is a unit for $\kappa$ precisely
because $\eta$ is a unit for $\mu$. The monoidality of $\delta_X$ and
$\varepsilon_X$ are due to the multiplicativity and unitalness of
$\Delta_H$ and $\varepsilon_H$. Finally, $G$ is a comonad because $H$
is a coalgebra.
\end{proof}

We now return to the question of dualizing Schauenburg's theorem. The
original proof relies heavily on the fact that a left bialgebroid $A$
is a generator in the category $_A\M$. This allows us to express the
coproduct of $A$ as the action on $1_A \oR 1_A$, $\Delta_A: a \mapsto
a_{(2)} \oR a_{(2)} := a \cdot (1_A \oR 1_A)$. An application of this
reasoning seems impossible. Consider, however, the following monadic
reformulation of the problem. A monoidal structure on $^H\M$ such that
$F: {^H\M} \to {^{C^e}\M} \simeq {\bicomC}$ is strict monoidal implies
that the monoidal product on ${^{C^e}\M}$ is \emph{lifted} to the
Eilenberg--Moore category of $\mathbb{G}$-coalgebras in the following
sense:
$$\xymatrix{
^{\mathbb{G}}\M \times {^{\mathbb{G}}\M} \ar[r]^-{\hat{\Co}} \ar[d]_-{F \times F}& ^{\mathbb{G}}\M \ar[d]^{F}\\
{^{C^e}\M} \times {^{C^e}\M} \ar[r]_-{\coC} & {^{C^e}\M}
}$$
This is a special case of the problem of 'liftings of functors',
orginally considered by Johnstone (\cite{Jon}). Our reference is
\cite{Wis} (this volume), from which we quote part (1) of Theorem 3.3.
\begin{thm}\label{robert}
Let $\mathbb{G} = \langle G, \delta, \varepsilon \rangle$ and  
$\mathbb{G'} = \langle G', \delta', \varepsilon' \rangle$ be comonads
on the categories $\M$ and $\M'$, respectively, and let $T: \M' \to \M$ be a functor. 
Denote $U: {^{\mathbb{G}}\M} \to \M$ and $U': {^{\mathbb{G}'}\M} \to \M'$ the canonical forgetful functors.

Then, the liftings $\hat{T}: {^{\mathbb{G}'}\M} \to {^{\mathbb{G}}\M}$ of $T$, in the sense:
$$\xymatrix{
{^{\mathbb{G}'}\M} \ar[r]^-{\hat{T}} \ar[d]_-{U'}& {^{\mathbb{G}}\M} \ar[d]^-{U}\\
\M' \ar[r]_-{T} & \M
}$$
are in bijective correspondance with natural
transformations $\kappa: T G' \to G T$ for which the following
diagrams commute:
$$\xymatrix{
T G' \ar[r]^-{T \delta'} \ar[d]_-{\kappa}& T G' G' \ar[r]^-{\kappa G} & G T G'
\ar[d]^-{G \kappa} & & T G' \ar[r]^-{T \varepsilon'} \ar[d]_-{\kappa}& T\\
G T \ar[rr]_-{\delta T} & & G G T & & G T \ar[ur]_-{\varepsilon T} &  
}$$
\end{thm}

Taking $\M' = {^{C^e}\M} \times {^{C^e}\M}$,  $\M = {^{C^e}\M}$ and 
$T = \_\, \coC \,\_\, :  {^{C^e}\M} \times {^{C^e}\M} \to {^{C^e}\M}$, 
we find that liftings of the monoidal structure to $^{\mathbb{G}}\M \simeq {^H\M}$
are in bijective correspondance with natural transformations
\begin{align}
\kappa_{M,N}: (H \coCe M) \coC (H \coCe N) &\to H \coCe (M \coC N) \nonumber
\end{align}
inducing commutative diagrams
$$\xymatrix{
G(M) \coC G(N) \ar[r]^-{\delta_M \Co \delta_N} \ar[d]_-{\kappa_{M,N}}& G^2(M) \coC G^2(N) 
\ar[r]^-{\kappa_{G(M),G(N)}} & G (G(M) \coC G(N)) \ar[d]^-{G \kappa_{M,N}} \\
G (M \coC N) \ar[rr]_{\delta_{M,N}} & & G^2 (M \coC N) 
}$$ 
\beq
\delta_{M\Co N} \ci \kappa_{M,N} = G \kappa_{M,N} \ci \kappa_{G(M),G(N)} \ci (\delta_M \coC \delta_N)
\eeq

and

$$\xymatrix{
G(M) \coC G(N) \ar[r]^-{\varepsilon \coC \varepsilon} \ar[d]_-{\kappa_{M,N}}& M \coC N \\
G (M \coC N) \ar[ur]_-{\varepsilon_{M\Co N}} &
}$$
\beq
\varepsilon_{M\Co N} \ci \kappa_{M,N} = \varepsilon_M \coC \varepsilon_N
\eeq
The two diagrams above recover two of the compatibility relations (\ref{a} and \ref{c}) of a bicomonad.
If, furthermore, we have an arrow $\xi: C \to G (C)$ making $C$ a $\mathbb{G}$--coalgebra such that the remaining
two bicomonad conditions (\ref{b} and \ref{d}) are satisfied, then $^{\mathbb{G}}\M$ becomes a (unital) monoidal 
category. Summarizing, we have the following weakened form of Schauenburg's theorem:
\begin{thm}
Let $\langle H, \bar{\Delta}, \bar{\varepsilon} \rangle$ be a comonoid in ${^{C^e}\M}$. Then there is a bijective
correspondance between
\begin{enumerate}
\item monoidal structures on ${^H\M}$ such that the forgetful functor $F: {^H\M} \to {^{C^e}}\M$ is strict monoidal 
\item a map $\kappa_{M,N}: (H \coCe M) \coC (H \coCe N) \to H \coCe (M \coC N)$, natural in both arguments and a 
map $\xi: C \to H \coCe C$ such that $\langle H, \bar{\Delta}, \bar{\varepsilon}; \kappa, \xi \rangle$ constitutes
a bicomonad, i.e.~the compatibilty conditions \ref{a}, \ref{b}, \ref{c} and \ref{d} are satisfied.   
\end{enumerate} 
\end{thm}  
Notice that in proving Theorem \ref{Hcomod}, we established (1) by constructing the maps $\kappa$ and $\xi$ of 
(2) from bicoalgebroid structure maps. For bialgebroids, a stronger result can be proved because a 
bialgebroid structure not only implies, but is equivalent to, the analogue of (2).     
    
\medskip
We now turn to the scalar extension of bicoalgebroids to investigate it from a comonadic
point of view. A scalar extension $H' = D \shrp H$ of the bicoalgebroid $H$ by the BCC $D$ gives
rise to an adjunction between the respective comodule categories. The forgetful functor
$F': {^{D\shrp H}\M} \to {^H\M}$ is induced by the epi $(\varepsilon_D \o H): D \shrp H \to H$.
It has the right adjoint induction functor
\begin{align}
& I': {^H\M} \to {^{D\shrp H}\M}, \;\; X \mapsto (D\shrp H) \coH X \nonumber
\end{align}
with the unit and counit of the adjunction being
\begin{align}
&\upsilon : X \to IF (X) = {D\shrp H} \coH X, \;\; x \mapsto \delta_{D\shrp H} (x) \nonumber \\
&\tau : FI (Y) = {^H(D\shrp H)} \coH Y \to Y, \;\; (d \o h) \o y \mapsto \varepsilon_D(d) \varepsilon(h) y
\nonumber
\end{align}

\medskip

As $C$--comodules, $D \shrp H = D \coC H$, so the induction functor
${D\shrp H} \coH \_$ is isomorphic to $D \coC \_\,$. By Lemma
\ref{BCClem}, this functor is strong monoidal, hence also
opmonoidal. It will remain opmonoidal upon composition with the
opmonoidal forgetful functor, making the canonical comonad $
\mathbb{G} = \langle F'I', F' \upsilon I', \tau \rangle$ an
opmonoidal endofunctor. The compatibility of the opmonoidal and
comonadic structure make $\mathbb{G}$ an opmonoidal comonad, not to
be confused with the monoidal comonad which we have christened
'bicomonad' earlier. We state the definition as concisely as
possible.

\begin{defi}
An opmonoidal comonad $\langle \langle G, G_{X,Y}, G^0 \rangle, \Delta, \varepsilon \rangle$ on a
monoidal category $\langle \mathcal{C}, \Co, \iota \rangle$ consists of

\begin{itemize}
\item An opmonoidal endofunctor $\langle G, G_{X,Y}, G^0 \rangle$ on $\M$ and
\item a comonad $\langle G, \Delta, \varepsilon \rangle$
\end{itemize}
such that $\Delta$ and $\varepsilon$ are opmonoidal
\end{defi}

\begin{pro}
Let $D$ be a BCC over the left bicoalgebroid $H$. Then the endofunctor $G = D \coC \_\,$ is an opmonoidal comonad on
$^H\M$.
\end{pro}
\begin{proof}
Recall that $G_{X,Y}: D \Co (X \Co Y) \to (D \Co X) \Co (D \Co Y)$ reads, on elements:
$d \o x \o y \mapsto d_{(1)} \o x\coa{-1} \o d_{(2)}\ract x\coa{-1} \o y$ and $G^0 = (\varepsilon_D \Co C):
D \Co C \to C$. The comonad structure follows from the coalgebra structure of $D$.

We have only to check the compatibility of the comonad and opmonoidal structure, meaning four commutative diagrams.
Opmonoidality of the comultiplication means (1):
$$\xymatrix{
D \Co (X \Co Y) \ar[r]^-{G_{X,Y}} \ar[d]_-{\Delta_{X \Co Y}}& (D \Co X) \Co (D \Co Y) \ar[d]^-{\Delta_X \Co \Delta_Y}\\
(D \Co D) \Co (X \Co Y) \ar[r]_-{G_{D \Co X, D \Co Y}\ci (D \Co G_{X,Y})}& (D \Co D \Co X) \Co (D \Co D \Co X)
}$$

An easy calculation shows the commutativity of the diagram. The upper and right hand side map compose to give

\begin{align}
d &\o x \o y \mapsto d_{(1)} \o (d_{(2)(1)} \o x\coa{0}) \o (d_{(2)(2)} \ract x\coa{-1} \o y) \mapsto
\nonumber\\
&\mapsto d_{(1)(1)} \o (d_{(2)(1)} \o x\coa{0})\coa{0} \o d_{(1)(2)} \ract (d_{(2)(1)} \o x\coa{0})\coa{-1} \o
d_{(2)(2)} \ract x\coa{-1} \o y \nonumber
\end{align}
\smallskip\\
Using the braided cocommutativity of $D$, we have:
\smallskip
\begin{align}\nonumber
&d_{(1)(1)} \o {d_{(2)(1)}}\coa{0} \o x^{\langle 0 \rangle\langle 0 \rangle}\o d_{(1)(2)}\ract 
({d_{(2)(1)}}\coa{-1} x^{\langle 0 \rangle\langle -1 \rangle}) \o d_{(2)(2)}\ract x\coa{-1} \o y = \\
&= d_{(1)(1)} \o d_{(1)(2)} \o x^{\langle 0 \rangle\langle 0 \rangle} \o
d_{(2)(1)}\ract x^{\langle 0 \rangle\langle -1 \rangle} \o d_{(2)(2)} \ract x\coa{-1} \o y, \nonumber
\end{align}
which is the composition of the lower and left hand side maps. The second diagram (2) for the opmonoidality of
$\Delta$ is
$$\xymatrix{
D \Co C \ar[dr]_-{G^0 \Co C} \ar[rr]^-{\Delta_D}& & D \Co D \Co C \ar[dl]^-{G^0 \ci
(D \Co G^0)}\\
      & C &
}$$
which commutes by the counit property of $\varepsilon_D$ ($G^0 = \varepsilon_D \Co C$). The remaining two diagrams
stating the opmonoidality of $\varepsilon_X : D \Co X \to X$ are (3)
$$ (\varepsilon_X \Co \varepsilon_Y) \ci G_{X,Y} = \varepsilon_{X,Y},$$
commuting since $\varepsilon(d_{(1)})x\coa{0} \coC \varepsilon(d_{(2)} \ract x\coa{-1}) y = \varepsilon(d)x \coC y$,
and (4)
$$ \varepsilon_C = G^0$$
which is a triviality.
\end{proof}

We briefly recall the dual situation, for scalar extensions of bialgebroids (for details, see \cite{Kor-Bim}).
For a scalar extension $Q \# H$ of a (left) bialgebroid H over $R$, the inclusion $\iota: H \into Q \# H$ induces
a monoidal forgetful functor $U: {_{Q\# H}\M} \to {_H\M}$. The \emph{left} adjoint of $U$ is
\begin{align}
& I: {_H\M} \to {_{Q\,\#\, H}\M}, \;\; X \mapsto (Q\,\#\, H) \otimes_H X \nonumber
\end{align}
with unit and counit
\begin{align}
&\eta : X \to UI (X) = {_H(Q\,\#\, H)} \otimes_H X, \;\; x \mapsto (1_Q \o 1_H)\o x \nonumber \\
&\varepsilon : IU (Y) = {_{Q\,\#\, H}(Q\,\#\, H)} \otimes_H Y \to Y, \;\; (q \o h) \o y \mapsto
(q \,\#\, h) \lact y = q \cdot (h\lact y)
\nonumber
\end{align}

Note that $U = \Hom_{H-} (H,\_\,)$. By Prop.~4.7 of \cite{Kor-Bim}, $I$ is strong monoidal, so the
underlying endofunctor of the canonical monad $\mathbb{T} = \langle
UI, U\varepsilon I, \eta\rangle$ on $_H\M$
will be monoidal, being the composition of two monoidal functors. Thus, the scalar extension of bialgebroids
gives rise to a monoidal monad on the module category of the 'smaller' bialgebroid.

\medskip
Having seen that scalar extensions of bialgebroids and
bicoalgebroids by BCA's and BCC's give rise to monoidal monads and
opmonoidal comonads, respectively, in the rest of this paper, we
make some tentative steps in the reverse direction.
\medskip

First, note that any monoidal category $\langle \mathcal{C}, \Co , I \rangle $ may be embedded
(monoidally, but not fully) into the category of it's endofunctors (\cite{Joyal-Street1}), which is 
monoidal with the composition of functors
as monoidal product and the identity functor as monoidal unit. The inclusion is given by $\mathcal{C} \into
\mathcal{C}^{\mathcal{C}}$, $X \mapsto \hat{X} \Co \_\;$, and the image of the inclusion will be denoted
$\hat{\mathcal{C}}$. The arrows of $\hat{\mathcal{C}}$ are natural
transformations of the form $\alpha_Z = \alpha \Co Z: X \Co Z \to Y
\Co Z$, with $\alpha : X \to Y$ an arrow in $\mathcal{C}$. An immediate
consequence is that for any map $\gamma \in \mathcal{C}$:
\beq
\gamma_{X\Co Y} = \gamma_X \Co Y,
\eeq
since $\gamma_{X\Co Y} = (\gamma \Co X) \Co Y = \gamma_X \Co Y$.

\begin{pro}
Let $\mathcal{C}$ be a monoidal category, and $\langle D\Co\_\,, D_{X,Y}, D^0\rangle$ an opmonoidal
endofunctor in $\hat{\mathcal{C}}$. Then there is a natural
transformation (between endofunctors of $\hat{\mathcal{C}}$),
$\hat{\theta}_Y: D \Co Y \Co \_\, \to Y
\Co D \Co \_\,$ such that $\langle G, \hat{\theta} \rangle$ is an object of $\lZ(\hat{\mathcal{C}})$.
\end{pro}
\begin{proof}
We shall only prove the latter statement, which amounts to constructing a natural transformation
\beq
\hat{\theta}_Y (\_):  (D \Co Y) \Co \_\, \to (Y \Co D) \Co \_
\eeq
satisfying eqs.~\ref{centrum1} and \ref{centrum2}. It is easily verified that

\begin{align}
\hat{\theta}_Y (X): (D \Co Y) \Co X &\to (D \Co Y) \Co  (D \Co X) \to (Y \Co D) \Co X \\
\hat{\theta}_Y (X)&= (G^0 \Co 1) \ci G_{Y,X} \nonumber
\end{align}
is appropriate. Eq.~\ref{centrum1} means that the following diagram
commutes (suppressing natural isomorphisms):\medskip

\xymatrix{
D \Co (X \Co Y) \Co Z \ar[r]^-{D_{X,Y} \Co Z} \ar[dd]_-{D_{X\Co Y, W}}& (D \Co X) \Co (D \Co Y)
\Co W \ar[r]^-{D^0 \Co 1} & X \Co D \Co (Y \Co W) \ar[d]^-{X \Co
  D_{Y,W}}\\
 &  & X \Co (D \Co Y) \Co (D \Co W) \ar[d]^-{X \Co D^0 \Co 1}\\
(D \Co X \Co Y) \Co (D \Co W) \ar[rr]_-{D^0 \Co 1} & & (X \Co Y) \Co D
\Co W
}

\bigskip

The diagram commutes by the coassociativity of ${D}_{X,Y}$ and because
$D^0$ is a counit for $D_{X,Y}$. Note that because $D_{X,Y} =
D_{X,\iota \Co Y} = D_{X,\iota} \Co Y$, the counit relation for $D^0$ is equivalent to the 
following property for $D_{X,Y}$:

$$\xymatrix{
D \Co X \Co Y \ar[r]^-{D_{X,Y}} \ar[d]_{D^0}& (D \Co X) \Co (D \Co Y) \ar[d]^-{(D^0
  \Co X) \Co (D^0 \Co Y)}\\
\iota \Co X \Co Y \ar[r] & (\iota \Co X) \Co (\iota \Co Y)
}$$

\end{proof}

\begin{rmk}
Clearly, by the embedding $\mathcal{C} \into \hat{\mathcal{C}}$, we
have in fact proven that $D$ is an object in $\lZ(\mathcal{C})$.
\end{rmk}
\begin{rmk}
Entirely analogously, for a monoidal category $\langle \M, \o, \iota \rangle$, and a \emph{monoidal} endofunctor
$\langle T, T_{X,Y}, T_0 \rangle$ in $\hat{\M}$, of the form $T = Q \o \_\,$ the natural transformation

\begin{align}
\hat{\bar{\theta}}_Y (X): Y \o Q \o X &\to Q \o Y \o Q \o X \to Q \o Y \o X\\
\hat{\bar{\theta}}_Y (X) &= T_{X,Y} \ci (T_0 \o 1) \nonumber
\end{align}
makes $\langle T, \hat{\bar{\theta}}_Y \rangle$ an object in $\rZ(\M)$
\end{rmk}

We saw that for a BCC $D \in {^H\YD_H}$ over a bicoalgebroid, $D
\coC \_\,$ is not only an opmonoidal endofunctor, but an opmonoidal
comonad. Unfortunately, it seems unlikely that the correspondence
between opmonoidal endofunctors of $\hat{\mathcal{C}}$ and objects
of $\lZ(\mathcal{C})$ can be extended to a correspondence between
opmonoidal \emph{comonads} of $\hat{\mathcal{C}}$ and  BCC's in
$\lZ(\mathcal{C})$ without further assumptions.

\section*{Acknowledgments}

I am indebted to Korn\'el Szlach\'anyi for his continued interest and 
our discussions, especially concerning the cocenter construction of \ref{cocentrum}.    
I'm thankful also to Gabi B\"ohm for her comments and for calling refs.~\cite{Jon} and \cite{Wis} to my 
attention. The suggestions of the Referee helped much to clarify the exposition.


\begin{thebibliography}{XX}
\begin{small}

\bibitem{Kor-Bim} I. B\'alint, K. Szlach\'anyi, \textit{Finitary
    Galois extensions over noncommutative bases}, J. Algebra \textbf{296} (2006)
520-560

\bibitem{Brz-Mil} T. Brzezi\'nski, G. Militaru, \textit{Bialgebroids,
$\times_A$-bialgebras and duality}, J. Algebra \textbf{251} (2002) 279-294

\bibitem{Brz-Wis} T. Brzezi\'nski, R. Wisbauer, \textit{Corings and Comodules},
London Math. Soc. LNS 309, Cambridge Univ. Press 2003

\bibitem{Cae-Wang-Yin} S. Caenepeel, Dingguo Wang, Yanmin Yin, \textit{Yetter-Drinfeld modules over weak 
Hopf algebras and the center construction}, Ann. Univ. Ferrara - Sez. VII - Sc. Mat. \textbf{51} (2005), 
69-98.  

\bibitem{Joyal-Street1} A. Joyal, R.H. Street, Braided tensor categories, \textit{Advances in Mathematics} 
\textbf{102}, 20-78 (1993)

\bibitem{Joyal-Street2} A. Joyal, R.H. Street, Tortile Yang-Baxter operators in tensor categories, 
\textit{J. Pure Appl. Algebra} \textbf{71}, 43-51 (1991)

\bibitem{Jon} Johnstone, P.T., \textit{Adjoint lifting theorems for categories of modules}, Bull. Lond. Mat.
Soc. 7, 294-297 (1975)

\bibitem{K-Sz} L. Kadison, K. Szlach\'anyi, \textit{Bialgebroid actions on depth
two extensions and duality}, Advances in Mathematics \textbf{179} (2003) 75-121

\bibitem{Lar-D2} L. Kadison, \textit{Co--depth two and related topics}, \texttt{math.QA/0601001}

\bibitem{Kassel} C. Kassel, ``Quantum Groups'', \textit{Grad. Textsin Math.} \textbf{155}, Springer Verlag, 
Berlin, 1995

\bibitem{MacLane} S. Mac Lane, \textit{Categories for the Working
Mathematician}, 2nd edition, GTM 5, Springer-Verlag New-York Inc., 1998

\bibitem{Majid} S. Majid, Representations, duals and quantum doubles of monoidal categories, 
\textit{Rend. Circ. Mat. Palermo (2)} Suppl. No. 26 (1991), 197-206   

\bibitem{Moer} I. Moerdijk, \textit{Monads on tensor categories}, Journal of Pure and Applied Algebra, 
\textbf{168} (2002) 189-208

\bibitem{SM:book} S. Montgomery, \textit{Hopf Algebras and Their Actions on Rings}, 
CBMS Lecture Notes \textbf{82}, Am. Math. Soc., Providence, RI (1993)

\bibitem{Schauenburg: ddqg} P. Schauenburg, \textit{Duals and
doubles of quantum groupoids} in "New trends in Hopf algebra theory",
Contemporary Mathematics {\bf 267},  p. 273, AMS 2000

\bibitem{Schn1} Hans--J\"urgen Schneider, \textit{Principal homogeneous spaces for
arbitrary Hopf algebras}, Isr. J. of Mathematics, Vol. 72, Nos. 1-2. (1990)

\bibitem{Sz:EilMoor} K. Szlach\'anyi, \textit{The monoidal Eilenberg--Moore construction and
bialgebroids}, Journal of Pure and Applied Algebra \textbf{182} (2003) 287--315

\bibitem{Tak:cohom} M. Takeuchi, \textit{Morita theorems for categories of comodules}, 
 J. Fac. Sci. Univ. Tokyo, Sec. IA, 24, (1977) 629-644 

\bibitem{CRTvO:cocenter} B. Torrecillas, F. van Oystaeyen, Y. H. Zhang,  
\textit{The Brauer Group of a Cocommutative Coalgebra}, Journal of  Algebra, \textbf{177} (1995) 536-568

\bibitem{Wis} Robert Wisbauer, \textit{Algebras versus coalgebras}, to appear in Applied Categorical 
Structures 
(this volume)

\end{small}
\end{thebibliography}
\end{document}